\input amstex.tex
\documentstyle{amsppt}
\overfullrule=0pt

\leftheadtext{C.D.Hill and M.Nacinovich}


\magnification=\magstep1

\font\sc=cmcsc10

\newcount\z
\newcount\p
\newcount\q
\newcount\r
\newcount\s
\newcount\t
\newcount\u

\def\thn{\number\p}
\def\sn{\number\q}
\def\tn{\number\r}
\def\lemn{\number\p}
\def\exn{\number\t}
\def\osn{\number\u}

\long\def\se#1{\advance\q by 1
\r=0  \p=0 \s=0 \t=0 \u=0
\bigskip
\noindent
\S\sn \quad
{\bf {#1}}\par
\nopagebreak}
\def\res{\roman{res}\,}
\def\singsupp{\roman{singsupp}\,}
\def\supp{\roman{supp}\,}
\def\ho{\roman{H}\,}
\def\Hom{\roman{Hom}}
\def\Res{\roman{Res}}
\long\def\thm#1{\advance\p by 1
\bigskip\noindent%
{\sc Theorem \sn .\thn}%
\quad{\sl #1} \par\noindent}

\long\def\prop#1{\advance\p by 1
\bigskip\noindent%
{\sc Proposition \sn .\thn}%
\quad{\sl #1} \par\noindent}

\long\def\thml#1#2{\advance\p by 1
\bigskip\noindent
{\sc Theorem \sn.\thn\, ({#1})}
\quad{\sl #2} \par\noindent}

\long\def\lem#1{\advance\p by 1
\medskip\noindent
{\sc Lemma \sn .\lemn}
\quad{\sl #1} \par\noindent}

\long\def\ex{\advance\t by 1
\medskip
\noindent
{\sc Example} \sn .\exn \quad}

\long\def\os{\medskip\advance\u by 1
\noindent
{\sc Remark} \sn.\osn \quad}

\long\def\form#1{\global\advance\r by 1
$${#1} \tag \sn.\tn$$}

\long\def\cor#1{\advance\p by 1
\bigskip\noindent%
{\sc Corollary \sn .\thn}%
\quad{\sl #1} \smallskip\noindent}
\def\dimo{\par\noindent {\sc Proof}\quad}

\topmatter
\title Leray Residues and Abel's Theorem in $CR$ codimension $k$
\endtitle
\author C.Denson Hill and Mauro Nacinovich
\endauthor
\address C.D.Hill - Department of Mathematics, SUNY at Stony Brook,
Stony Brook NY 11794, USA \endaddress
\address M.Nacinovich - Dipartimento di Matematica "L.Tonelli" -
via F.Buonarroti, 2 - 56127 PISA, Italy \endaddress
\email Dhill\@ math.sunysb.edu \quad nacinovi\@ dm.unipi.it\endemail
\subjclass 35 32 53 \endsubjclass
\endtopmatter
\document
\centerline{\sc Contents}\par
\roster
\item"\S 1"  Introduction
\item"\S 2"  Preliminaries and notation
\item"\S 3"  Local defining functions
\item"\S 4"  The residue form
\item"\S 5"  Properties of the residue form
\item"\S 6"  The residue formula
\item"\S 7"  The homological residue
\item"\S 8"  The cohomological residue
\item"\S 9"  Properties of the residue class and global defining
functions
\item"\S 10" Higher order poles
\item"\S 11" The calculus of residues
\item"\S 12" Iterated residues
\item"\S 13" The calculus of residues for intersecting poles
\item"\S 14" Abel's global residues theorem
\item"\S 15" Applications of the Abel theorem
\endroster

\se{Introduction}
In this paper we generalize Leray's calculus of residues in several
complex variables [Lr], to the situation of an {\it abstract}
smooth $CR$ manifold $M$ of general type $(n,k)$, and a {\it polar}
submanifold $S$ which is a {\it smooth} $CR$ submanifold of type
$(n-1,k)$, transversal to the Levi distribution of $M$. Here
$n$ is the $CR$ dimension of $M$ and $k$ its $CR$ codimension; so
$\roman{dim}_{\Bbb R}M=2n+k$. When $k=0$ this means that $M$
is an $n$-dimensional complex manifold, and $S$ is a complex
submanifold having complex codimension $1$ in $M$; for $k=0$
and $n=1$, everything reduces to the classical theory of residues
in one complex variable.\par
The extension of the Leray residue calculus to abstract $CR$ manifolds
$M$ of $CR$ codimension $k$ involves a certain number of new
turns and twists, which require investigation. One of these is
the strange fact that the polar submanifold $S$ may not have local
defining functions that are $CR$ on $M$ (see the example at the
end of \S 3). In particular this forces upon us an enlargement
of the notion of semimeromorphic function or form, allowing
singularities along the polar submanifold $S$ that are
more general than
what one may expect by analogy with the complex manifold case.
Another new aspect is that when one takes a maximal number $n$
of polar submanifolds $S_1,\, S_2,\,\hdots,\, S_n$, having
normal crossings, the intersection $S=S_1\cap S_2\cap\cdots\cap S_n$
is of type $(0,k)$, and hence totally real. We then obtain
$(k+1)$ different kinds of "Grothendieck" residues, which
are the analogues of the Grothendieck point residue that one
has when $k=0$ (see [D]).
\par
In the case where each polar submanifold $S_j$ has a
global defining function in a neighborhood of $S$, we obtain
an actual {\it calculus of residues} which entails using only
the usual operations of the exterior differential calculus
of smooth forms (see \S 13). In this respect our residue
calculus improves somewhat that of Leray, even for the
complex manifold case ($k=0$).
\par
In \S 14 we use the theory we have developed to generalize the
classical theorem of Abel, about the sum of the residues of an
abelian differential of the second kind on a compact Riemann
surface, to the case of a compact {\it abstract} $CR$
manifold $M$. In this connection our discussion follows
along the lines of Griffiths [G]; however we actually obtain
a more general result, even for the complex manifold case ($k=0$),
because we allow for the intersection of only $m$ ($1\leq m\leq n$)
polar submanifolds having normal crossings. Thus we obtain
(Theorem 14.2) a result that applies to forms of more general
degree. In the last section we apply the Abel theorem to
derive quite general {\it period relations} (Propositions
15.1 and 15.2).
\se{Preliminaries and notation}
An {\it abstract $CR$ manifold} of type $(n,k)$ is a triple
\footnote{But we shall write only $M$ for simplicity when
no confusion can arise.}
$(M,HM,J)$
where $M$ is a paracompact smooth manifold of
dimension $2n+k$, $HM$ is a smooth subbundle of even dimension
$2n$ of the tangent bundle $TM$ of $M$, which is called
the {\it Levi distribution},
and $J$ is a smooth complex
structure on the fibers of $HM$: this means that $J:HM @>>> HM$
is an equivalence of smooth vector bundles with $J^2=-Id$.
We also require that $J$ be {\it formally integrable}. This
condition can be expressed in terms of the complex subbundle
\form{T^{0,1}M\, =\, \left\{\left. X\,
+\,{\ssize{\sqrt{-1}}}\,JX\,\right|
X\in HM\right\}}
of the complexified tangent bundle $\Bbb CTM$ of $M$, by requiring that
\form{\left[ \Gamma(M,T^{0,1}M), \Gamma(M,T^{0,1}M)\right]
\subset \Gamma(M,T^{0,1}M)\, .}
We note that $T^{0,1}M$ is the eigenspace corresponding to the
eigenvalue $-{\ssize{\sqrt{-1}}}$ of $J$.
Its complex conjugate with
respect to the real form $TM$ of $\Bbb CTM$:
\form{T^{1,0}M\,=\,
\left\{\left. X\, - \,
{\ssize{\sqrt{-1}}} JX\, \right| \, X\in HM\right\}}
is the eigenspace corresponding to the eigenvalue
${\ssize{\sqrt{-1}}}$ of $J$.
We have therefore
\form{T^{1,0}M=\overline{T^{0,1}M}, \qquad
T^{1,0}M\cap T^{0,1}M\, = \, 0_M\;(\text{zero section of $\Bbb CTM$})\, .}
We note that the datum of a complex subbundle $T^{0,1}M$
of rank $n$ of the
complexified tangent bundle $\Bbb CTM$, satisfying (2.2) and (2.4),
defines on the smooth manifold $M$ of dimension $2n+k$ a unique
structure of $CR$ manifold of type $(n,k)$: indeed we define
\form{HM\, = \,
\left\{\left. \Re Z\, \right|\, Z\in T^{0,1}M\right\}}
and, for $X\in HM$, we set $JX=Y$ iff
$X+{\ssize{\sqrt{-1}}}Y\in T^{0,1}M$.
The map $J:HM@>>>HM$ is well defined because of (2.4).
\par
An abstract $CR$ manifold of type $(n,0)$ is a complex manifold
of dimension $n$ by the Newlander-Nirenberg theorem.
\par
Note that the vector bundle $HM$ can be considered a
{\it complex} vector bundle of rank $n$ on $M$, with
the complex structure on the fibers defined by $J$.\par
The annihilator bundle $H^0M$ of $HM$ in the cotangent bundle $T^*M$
of $M$ is called the {\it characteristic bundle} of $M$. The
quotient bundle $\left. T^*M\right/ H^0M$ is the dual bundle
of $HM$. Therefore it is a complex vector bundle of rank $n$
and hence orientable, with the natural orientation associated
to its complex structure. It follows that {\sl the abstract
$CR$ manifold $M$ is orientable iff its characteristic
bundle $H^0M$ is orientable}.
\smallskip
We use the notation $\Cal E(M)$ for the {\it exterior algebra}
of smooth complex valued alternating forms on $M$:
\form{\Cal E(M)\, = \, \bigoplus_{p=0}^{2n+k}{\Cal E^p(M)}\, ,}
where $\Cal E^p(M)$ is the subspace of forms homogeneous
of degree $p$. In particular $\Cal E^0(M)=\Cal C^\infty(M)=
\Cal C^\infty(M,\Bbb C)$.
\par
We consider the {\it ideal} in $\Cal E(M)$:
\form{\Cal I(M)\, =\, \left\{
\alpha\in \Cal E(M)\, \left| \, \alpha\left|_{T^{0,1}M}\,=\,0
\right. \right.
\right\}\subset \bigoplus_{p=1}^{2n+k}{\Cal E^p(M)}\, }
and its exterior powers:
\form{\Cal I^0(M)=\Cal E(M),\quad \Cal I^p(M)=\Cal I^{p-1}(M)\wedge
\Cal I(M)\quad (p=1,2,\hdots)\, .}
The integrability condition (2.2) can also be expressed in terms
of the ideal $\Cal I(M)$ by:
\form{d\Cal I(M)\, \subset\, \Cal I(M)\, .}
We also have
\form{d\Cal I^p(M)\, \subset\, \Cal I^p(M)\quad \forall p=1,2,\hdots}
We note that $\Cal I^{n+k+1}=0$ by reasons of degree. Hence we have
a decreasing sequence of closed ideals:
\form{\Cal E(M)=\Cal I^0(M) \supset \Cal I^1(M) \supset
\Cal I^2(M)\supset \cdots \supset \Cal I^{n+k}(M) \supset \{0\}\, .}
For integers $0\leq p\leq n+k$, $0\leq q \leq n$ we set
\form{\Cal I^{(p,q)}M\, = \, \Cal I^p(M)\cap\Cal E^{p+q}(M)\, }
Then clearly we have:
\form{d\Cal I^{(p,q)}(M)\subset\Cal I^{(p,q+1)}(M), \qquad
\Cal I^{(p+1,q)}(M)\subset\Cal I^{(p,q-1)}(M)\, .}
Now we define, for $0\leq p\leq n+k$, $0\leq q\leq n$:
\form{Q^{p,q}(M) \, = \,
\Cal I^{(p,q)}(M)\left/ \Cal I^{(p+1,q-1)}(M)\right.}
so that, passing to the quotient from the de Rham complex
we obtain the {\it Cauchy-Riemann complexes}:
\form{\CD
0\rightarrow Q^{p,0}(M)@>{\bar\partial_M}>>
Q^{p,1}(M)@>{\bar\partial_M}>>
\cdots
@>{\bar\partial_M}>>
Q^{p,n}(M)\rightarrow 0\, ,\endCD}
for each $0\leq p\leq n+k$.
\par
We denote by $\Omega_M^p(M)$ the kernel of the map:
$\bar\partial_M:Q^{p,0}(M) @>>> Q^{p,1}(M)$. Note that
$Q^{0,0}(M)=\Cal E^0(M)$, and in general
$Q^{p,0}(M)\subset\Cal E^p(M)$. Thus $\Omega_M^p(M)\subset
\Cal E^p(M)$ and its elements are {\sl exterior differential
forms} (with complex valued coefficients) in $M$. When
$p=0$ we write $\Cal O_M(M)$ for $\Omega^0_M(M)$.
The functions $f$ in $\Cal O_M(M)$ are called $CR$ functions
on $M$ and the forms in $\Omega^p_M(M)$ $CR$ forms of degree
$p$ in $M$. \par
We note that any nonempty
open subset $U$ of $M$ is in a natural way
a $CR$ manifold of the same $CR$ dimension and $CR$ codimension.
We can therefore consider in a consistent way the space
$\Cal O_M(U)$ of $CR$ functions on $U$, the spaces
$\Omega_M^p(U)$ of $CR$ $p$-forms on $U$, etc.
\par
If $S$ is a {\sl real} submanifold of $M$ and $\psi$ a differential
form defined on a neighborhood of $S$, we denote by
$\psi\left|_S\right.$ the {\it pullback} of $\psi$ to $S$.
\par
We say that $S\subset M$ is a $CR$ submanifold of $M$ if
$$HS=\left(TS\cap HM\right)\cap J\left(TS\cap HM\right)$$
is a distribution of constant rank in $TS$. The triple
$\left(S,HS,J\left|_{HS}\right.\right)$ satisfies indeed
in this case the requirements for an {\sl abstract} $CR$
manifold of the definition at the beginning of the section.
\smallskip
Let $(M_1,HM_1,J_1)$ and $(M_2,HM_2,J_2)$ be two (abstract)
$CR$ manifolds of type $(n_1,k_1)$ and $(n_2,k_2)$ respectively.
A differentiable map $F:M_1@>>>M_2$ is called a $CR$ map
iff: ($i$) $dF(HM_1)\subset HM_2$, and ($ii$) $dF(J_1X)=
J_2dF(X)$ for every $X\in HM_1$.
\par

\se{Local defining functions}
In the sequel we shall consider the following situation:
$M$ is a smooth ($\Cal C^\infty$) {\it abstract} $CR$ manifold
of type $(n,k)$, where $n$ is the $CR$ dimension and $k$ the
$CR$ codimension. We assume that $M$ is connected, paracompact
(countable at infinity) and orientable (see \S 2). \par
$S$ will be a {\it polar submanifold} in $M$: By polar we mean
that $S$ is a smooth closed $CR$ submanifold of type $(n-1,k)$,
and that $S$ is transversal to the Levi distribution on $M$; i.e.,
$$ T_xS+H_xM\, = \, T_xM\, , \quad \forall x\in S\, .$$
Here {\sl closed} means as a subset; hence the topology of
the differentiable submanifold $S$ agrees with the one induced
on $S$ from the topology of $M$, and $S$ is also paracompact.
These assumptions imply that the characteristic bundle
$H^0S$ of $S$ is just the restriction to $S$ of the characteristic
bundle $H^0M$ of $M$. Since $M$ being orientable means that
$H^0M$ is orientable, we have that $H^0S$ is orientable;
hence a given orientation of $M$ induces a corresponding orientation
of $S$.
\lem{Let $p\in S$. Then there is an open neighborhood $U$ of $p$
in $M$ and a $\Cal C^\infty$ function $s:U @>>> \Bbb C$
with $S\cap U\, = \, \left\{x\in U\left| s(s)=0\right.\right\}$
and $ds\neq 0$ in $U$ satisfying
\roster
\item"($i$)" $\bar\partial_Ms=0$ on $S\cap U$;
\item"($ii$)" in fact, $s$ can be chosen so that
\form{ds\, = \, s\gamma\, + \, \eta\qquad\text{in}\quad U\, ,}
where $\gamma,\, \eta\in\Cal E^1(U)$ and $\eta\in\Cal I(U)$.
Moreover, the class $[\gamma]\in Q^{0,1}(U)$ defined by
$\gamma$ satisfies $\bar\partial_M[\gamma]=0$ on $U$.
\endroster}
\dimo In a sufficiently small neighborhood $U$ of $p$
there are local $\Cal C^\infty$ real valued defining functions
$\rho_1,\, \rho_2$ for $S$ with $d\rho_1,\, d\rho_2$ linearly
independent at each point of $U$. But since $S$ is polar,
$d\rho_1$ and $d\rho_2$ are not linearly independent
modulo the ideal $\Cal I(M)$ along $S$. Thus near $p$ there
is a uniquely determined complex valued smooth function
$f(x)$, defined along $S$, such that
$d\rho_1(x)+f(x)\cdot d\rho_2(x)\in\Cal I_x(M)$ on $S$.
Take any smooth extension $\tilde f$ of $f$ to $U$. Then
$s=\rho_1+f\rho_2$ satisfies ($i$).\par
We denote the function $s$ obtained above by $s_0$; because
of ($i$) it satisfies
\form{ds_0\, = \, s_0\gamma_0+\bar s_0\alpha_0+\eta_0\, ,\quad
\eta_0\in \Cal I(U)\, .}
Here and also below $\gamma_k,\, \alpha_k,\, \eta_k\in\Cal E^1(U)$
and $\eta_k\in\Cal I(U)$, and et cetera with primes. Applying $d$
to (3.2), we find that $d\bar s_0\wedge\alpha_0$ belongs to
$\Cal I(M)$ along $S$. By Cartan's lemma
$$\alpha_0=g_0 d\bar s_0+s_0\alpha_0'+\bar s_0\alpha_0''+
\eta_0'\, ,$$
where $g_0\in\Cal E^0(U)$. We set $s_1=s_0-\frac{1}{2}g_0(\bar s_0)^2$
and obtain a new local defining function in $U$ which satisfies
\form{ds_1\, = \, s_1\gamma_1+(\bar s_1)^2\alpha_1+\eta_1\, ,
\quad \eta_1\in\Cal I(U)\, .}
By induction we obtain a sequence $\left\{ s_m\right\}_{m=0}^\infty$
of local defining functions in $U$ satisfying
\form{ds_m\, = \, s_m\gamma_m+(\bar s_m)^{m+1}+\eta_m\, ,\quad
\eta_m\in\Cal I(U)}
and
\form{s_m-s_{m-1}=O\left(|s_0|^{m+1}\right)\quad\text{as}\quad
s_0 @>>> 0\, .}
We have constructed $s_0, \, s_1$ satisfying (3.4), (3.5) for
$m=1$. So assume we have $s_0,\, s_1, \, \hdots,\, s_m$ and to
verify the induction step, we construct $s_{m+1}$: applying
$d$ to (3.4) we obtain
\form{s_m d\gamma_m+(\bar s_m)^m\left[
(m+1)d\bar s_m\wedge\alpha_m+\bar s_m\alpha_m\wedge\gamma_m+
\bar s_m d\alpha_m\right]+\eta_m\wedge\gamma_m+d\eta_m\, = \, 0\, ,}
in which the last two terms belong to $\Cal I(U)$. \par
In order to exploit (3.6) we note that by the local triviality
of vector bundles, we may shrink $U$ if necessary, and obtain
\form{\Cal E^2(U)\left/_{\dsize{\Cal I(U)\cap\Cal E^2(U)}}\right.
\simeq \left[\Cal E^0(U)\right]^{\binom{n}{2}}\, .}
We consider the equations
\form{\cases
F=-(\bar s_m)^m V\\
G= s_m V\, ,
\endcases}
which are to be solved for $V\in\left[\Cal E^0(U)\right]^{\binom{n}{2}}$,
given $F,G\in\left[\Cal E^0(U)\right]^{\binom{n}{2}}$ which
satisfy the compatibility condition
\form{s_m F\, + \, (\bar s_m)^m G\, = \, 0\, .}
We may choose a coordinate system in $U$ of the form
$\left(s_m,\bar s_m,x_3,\hdots,x_{2n+k}\right)=
\left(s_m,\bar s_m,x\right)$. We have the exact sequence:
\form{\CD
\Bbb C[s_m,\bar s_m,x]@>{\left[{\matrix
-(\bar s_m)^m\\
s_m
\endmatrix}
\right]}>>
\Bbb C[s_m,\bar s_m,x]^2 @>{\left[s_m,(\bar s_m)^m\right]}>>
\Bbb C[s_m,\bar s_m,x]
\endCD}
of homomorphisms over the ring of polynomials in the
coordinates. Since the ring of formal power series is
flat over the ring of polynomials, the equation (3.8)
has a formal power series solution at each point of $U$.
Since $s_m$ and $(\bar s_m)^m$ may be regarded as real analytic
(polynomial) functions of the coordinates $x_1=\Re s_m$,
$x_2=\Im s_m$, we may apply the Whitney theorem on closed
ideals [T]; to obtain that (3.8) admits a smooth solution
$V$ in $U$. Returning to (3.6) we take $F$, $G$ to be the
projections of $d\gamma_m$ and
$\left[(m+1)d\bar s_m\wedge\alpha_m+\bar s_m\alpha_m\wedge\gamma_m
+\bar s_m d\alpha_m\right]$ into the quotient (3.7),
respectively. Then by the discussion above, we have a solution
$V$ to (3.8), and it follows that $d\bar s_m\wedge\alpha_m$
belongs to the ideal $\Cal I(M)$ along $S$. Again by
Cartan's lemma
$$\alpha_m=g_m d\bar s_m+s_m\alpha'_m+\bar s_m\alpha''_m+
\eta'_m\, ,$$
with $g_m\in\Cal E^0(U)$ and $\eta'_m\in\Cal I(U)$. Finally
we set $s_{m+1}=s_m-\frac{1}{m+2}g_m(\bar s_m)^{m+2}$. This
$s_{m+1}$ satisfies (3.4) and (3.5).\par
In a neighborhood $U$ of $p$ in $M$ we may now construct
a $\Cal C^\infty$ function $\tilde s$ such that
\form{\tilde s - s_m = O\left(|s_0|^{m+1}\right)\quad\text{as}\quad
s_0 @>>>0\, .}
This just boils down to the well-known fact that it is possible
to construct a $\Cal C^\infty$ function whose normal derivatives
are all smoothly prescribed along the smooth manifold $S$.
By construction this $\tilde s$ satisfies (3.1) for suitable
$\gamma,\, \eta\in\Cal E^1(U)$ with $\eta\in\Cal I(U)$. \par
Thus $d\gamma$ is in the ideal off of $S$, and hence in the
ideal across $S$, as $\gamma$ is smooth. Therefore
$\bar\partial_M[\gamma]\,=\,0$ in $U$, and the proof is complete.
\medskip
\noindent
{\sc Remark 1}.\quad  Note that the form $\gamma$ in (3.1) is
determined modulo $\Cal I(U)$.
If there is a neighborhood $U$ of $p$ in $M$
in which we can solve the equation
\form{\bar\partial_Mu=[\gamma]\qquad\text{in}\quad U}
for a function $u\in\Cal E^0(U)$, then $\hat s=e^{-u}\tilde s$
is a local defining function for $S$ which satisfies
$\bar\partial_M\hat s=0$ in $U$, i.e. there is a local
defining function for $f$ at $p$ which is $CR$ on $M$ in
a neighborhood of $p$.
\medskip\noindent
{\sc Remark 2}.\quad In particular, Remark 1 applies if $M$ is locally
embeddable at $p$ and is $2$-pseudoconcave there; indeed the
$\bar\partial_M$-Poincar\'e lemma for $Q^{0,1}$-forms on $M$ is
valid at $p$. (see  [N1], [HN1], [HN2]).
\lem{Suppose our abstract $M$ is $1$-pseudoconcave at $p$.
Let $s$, $\tilde s$ be two local $\Cal C^\infty$ defining functions
for $S$ at $p$, each satisfying (3.1). Then there is a local
nonvanishing $\Cal C^\infty$ function $h$ at $p$ such that
$\tilde s\, = \, hs$.}
\dimo We shall use the hypoellipticity of the $\bar\partial_M$
operator on functions; namely, if $u\in L^2_{loc}$ near $p$,
and if $\bar\partial_Mu$ is smooth in a neighborhood of $p$,
then $u$ is $\Cal C^\infty$ at $p$. This result is a consequence
of the interior subelliptic estimate with loss of $\frac{1}{2}$
derivative proved in [HN1], although the hypoellipticity
consequence was not explicitly stated there.
\footnote{For a proof of the hypoellipticity under the weaker
notion of {\it essential pseudoconcavity} see [HN3]; actually
Lemma 3.2. remains valid under this weaker assumption.}
\par
Set $u=\log\left(\frac{\tilde s}{s}\right)$. Then $u\in L^2_{loc}$
and it follows from (3.1) that $\bar\partial_Mu$ is smooth;
indeed $\bar\partial_Mu=[\tilde\gamma]-[\gamma]$, where as above
$[\,\cdot\,]$ denotes the class in $Q^{0,1}$. Hence $u$ is smooth
and we may take $h$ to be $e^u$.
\medskip
We end this section with an example showing that {\sl polar
submanifolds $S$ in $M$ do not always have local defining
functions that are $CR$.}
\smallskip\noindent
{\sc Example}.\quad Consider $M=\Bbb R^5=\Bbb C_z\times\Bbb C_w
\times\Bbb R_t$. We specify on $M$ the structure of an abstract $CR$
manifold of type $(2,1)$ by prescribing the following basis for
$T^{0,1}M$:
\form{\cases
\dfrac{\partial}{\partial\bar z}\,-\,
{\ssize{\sqrt{-1}}}\, z\dfrac{\partial}{\partial
t}\,+\, w f(z,t)\dfrac{\partial}{\partial w}\, ,\\
\dfrac{\partial}{\partial\bar w}\, .
\endcases}
We choose the smooth complex valued function $f(z,t)$ such that
the equation of Hans Lewy [L],
\form{\dfrac{\partial u}{\partial\bar z}\, - \,
{\ssize{\sqrt{-1}}}\, z\dfrac{\partial u}{\partial t}\, = \, f\, ,
\qquad\text{in}\quad \Bbb R^3\, ,}
has no local solution in a neighborhood of any point.
Then $S=\{w=0\}$ is a polar submanifold in $M$. Suppose that
$M$ could be locally defined, near some point, by a smooth
defining function $s$. In particular, $s$ satisfies
\form{\dfrac{\partial s}{\partial\bar z}\, - \,
{\ssize{\sqrt{-1}}}\, z\,\dfrac{\partial s}{\partial t}\, + \,
w f\dfrac{\partial s}{\partial w}\, = \, 0\, .}
Applying $\partial\left/\partial w\right.$ to (3.15)
yields
\form{\left(\dfrac{\partial}{\partial\bar z}-
{\ssize{\sqrt{-1}}}\, z
\dfrac{\partial}{\partial t}\right)\dfrac{\partial s}{\partial w}
+ f\dfrac{\partial s}{\partial w}+wf\dfrac{\partial^2s}{\partial w^2}
\,=\, 0\, .}
Note that $\partial s\left/ {\partial w}\right.\neq 0$ along $S$,
as $s$ is a defining function. Consider the function
\form{u(z,t)\, = \, -\log\dfrac{\partial s}{\partial w}(z,0,t)\, ,}
using any local branch of the logarithm. It gives a smooth
solution of (3.14), which contradicts our choice of $f$.
Note that this example gives another interpretation of the
example in [H].

\se{The residue form}
We return to the situation with $M$ and $S$ as at the beginning of
\S 3. Consider a form $\phi\in\Cal E^p(M\setminus S)$.
\smallskip
\noindent
{\sc Definition}.\quad {\sl $\phi$ is said to have a pole of
(at most) the first order along $S$ iff: given any point $p\in S$,
there exists an open neighborhood $U$ of $p$ in $M$, and a smooth
local defining function $s$ for $S$ in $U$, with $ds\neq 0$ in $U$,
and satisfying $\bar\partial_Ms=0$ on $S\cap U$, such that
$s\phi\in\Cal E^p(U)$.}
\thm{Let $\phi\in\Cal E^p(M\setminus S)$ be closed on $M\setminus S$
and have a first order pole along $S$. Let $p$, $U$, $s$ be as in
the above definition. Then
\roster
\item There exists $\psi\in\Cal E^{p-1}(U)$ and $\theta\in\Cal E^p(U)$
such that
$$\phi\, = \, \dfrac{ds}{s}\wedge \psi\, + \, \theta\, \quad
\text{in}\quad U\, .$$
\item $\psi\left|_S\right.$ is closed on $S\cap U$ and depends only
on $\phi$ and $S$.
\item If $\phi\in\Omega^p_M(M\setminus S)$ then
$\psi\left|_S\right.\in\Omega^{p-1}_S(S\cap U)$.
\endroster}
\dimo Since $\phi$ is closed, $d(s\phi)=ds\wedge\phi+sd\phi=
ds\wedge\phi$, showing that $ds\wedge\phi$ has the smooth
extension $d(s\phi)$ across $S\cap U$. Hence by continuity
we have $ds\wedge d(s\phi)=0$ in $U$. By Cartan's lemma
there exists $\theta\in\Cal E^p(U)$ with $d(s\phi)=ds\wedge\theta$
in $U$. It follows that $ds\wedge\left(s\phi-s\theta\right)=0$ in $U$.
Applying Cartan's lemma again, we obtain a $\psi\in\Cal E^{p-1}(U)$
such that $s\phi-s\theta=ds\wedge\psi$ in $U$. This establishes
the first point in the theorem.\par
Next we show that $\psi\left|_S\right.$ depends only on $\phi$
and on $s$: Suppose we have
$$\phi=\dfrac{ds}{s}\wedge\psi_1+\theta_1=
\dfrac{ds}{s}\wedge\psi_2+\theta_2$$
with $\psi_1,\,\psi_2\in\Cal E^{p-1}(U)$ and
$\theta_1,\, \theta_2\in\Cal E^p(U)$. Setting
$\psi=\psi_1-\psi_2$, $\theta=\theta_1-\theta_2$ we
subtract and multiply by $s$ to obtain
$ds\wedge\psi+s\theta=0$. Hence $ds\wedge(s\theta)=0$
and therefore $ds\wedge\theta=0$ in $U$. Once again
this yields $\theta=ds\wedge\omega$ in $U$, where
$\omega\in\Cal E^{p-1}(U)$. Substituting back in for $\theta$
we find that $ds\wedge\left(\psi+s\omega\right)=0$ in $U$.
As $\psi+s\omega$ is smooth in $U$, $\psi+s\omega=ds\wedge
\tilde\omega$ for an $\tilde\omega\in\Cal E^{p-2}(U)$. Hence
$\psi\left|_S\right.=0$ on $S\cap U$, so ${\psi_1}\left|_S\right.
={\psi_2}\left|_S\right.$, as claimed.\par
Return to the $\psi$ of the theorem.
To show that $\psi\left|_S\right.$
is closed, we observe that
$0=d\phi=-\dfrac{ds}{s}\wedge d\psi+d\theta$.
This means we can apply the argument above (to $d\phi$) and conclude
that $d\left(\psi\left|_S\right.\right)=(d\psi)\left|_S\right.=0$
on $S\cap U$. \par
Finally we show that $\psi\left|_S\right.$ does not depend on
the choice of the local defining function $s$. Let $s^*$ be
another smooth local defining function for $S$ in $U$, with
$ds^*\neq 0$ in $U$ and $\bar\partial_Ms^*=0$ on $S\cap U$,
such that $s^*\phi\in\Cal E^p(U)$. Using the defining function
$s^*$ we obtain
\form{\phi=\dfrac{ds^*}{s^*}\wedge\psi^*+\theta^*\, ,
\quad \psi^*\in\Cal E^{p-1}(U)\, ,\quad
\theta^*\in\Cal E^p(U)}
and also write
\form{\phi=\dfrac{ds^*}{s^*}\wedge\psi+\theta+
\left[\dfrac{ds}{s}-\dfrac{ds^*}{s^*}\right]\wedge\psi\, .}
If $s^*=hs$ where $h\in\Cal C^\infty(U)$, $h\neq 0$, then
$\left[\dfrac{ds}{s}-\dfrac{ds^*}{s^*}\right]=-\dfrac{dh}{h}$,
which is smooth across $S$, so the uniqueness argument
given above shows that ${\psi^*}\left|_s\right.=\psi\left|_S\right.$.
In general we have that $ds^*=hds$ along $S$, with
$h\in\Cal C^\infty(U)$, $h\neq 0$ on $U$: indeed
$\left(\Bbb CT^*_SM\right)_x\cap\Cal I_x(M)\simeq\Bbb C$ for
$x\in S\cap U$, where $\Bbb CT^*_MS$ denotes the complexification
of the conormal bundle of $S$ in $M$. By the discussion above
we may replace $s$ by $hs$ without changing $\psi\left|_S\right.$.
With this done, we now have that $ds^*=ds$ along $S\cap U$. It
follows that $s^*-s=O\left(|s|^2\right)$, which in turn implies
that $\left[\dfrac{ds}{s}-\dfrac{ds^*}{s^*}\right]$ is bounded
in $U$. Subtracting (4.2) from (4.1) we obtain
$0=\dfrac{ds^*}{s^*}\wedge\left[\psi^*-\psi\right]+\tilde\theta$
with $\tilde\theta$ bounded in $U$. This implies that
$ds^*\wedge\left[\psi^*-\psi\right]=0$ along $S$. This yields
that $\psi^*-\psi=ds^*\wedge\omega+s^*\omega'
+\overline{s^*}\omega''$
in $U$, where $\omega\in\Cal E^{p-2}(U)$ and $\omega',\,
\omega'' \in \Cal E^{p-1}(U)$. Therefore ${\psi^*}\left|_S\right.
=\psi\left|_S\right.$, completing the proof of the second
part of the theorem.
\par
Finally we establish the third point of the theorem: Since
$\phi\in\Omega^p_M(M\setminus S)$, we have in particular that
$\phi\in\Cal I^p(M\setminus S)$. By continuity $s\phi\in
\Cal I^p(M)$. Thus $d(s\phi)\in\Cal I^p(M)$ and we can choose
$\theta\in\Cal I^{p-1}(U)$ in (1).
Therefore $ds\wedge\psi\in\Cal I^p(M)$ along
$S\cap U$. This forces $\psi$ to belong to $\Cal I^{p-1}(U)$
along $S\cap U$. Since $\Cal I^{p-1}(S)=\Cal I^{p-1}(M)\left|_S\right.$
we conclude that $\psi\left|_S\right.\in\Cal I^{p-1}(S\cap U)$.
But $d\psi\left|_S\right.=0$, so $\psi\left|_S\right.
\in\Omega^{p-1}_S(S\cap U)$. This completes the proof of the theorem.
\bigskip
Let $F$ be a closed subset of $M$ which satisfies
\form{S\cap F\, \subset\, \overline{\left(M
\setminus S\right)\cap F}\, .}
Suppose $\phi=0$ on $F$. Then it is clear from the above costruction
that $\psi\left|_S\right.$ vanishes on $\left(S\cap F\right)\cap U$.
\par
Let $\phi\in\Cal E^p(M\setminus S)$ be closed on $M\setminus S$
and have a first order pole along $S$. Let $U$, $\psi$ as
in Theorem 4.1.
\smallskip
\noindent
{\sc Definition}.\quad {\sl The {\bf residue form} of $\phi$
on $S$ is locally defined by
\form{\res[\phi]\, = \, \psi\left|_S\right. \quad
\text{on}\quad S\cap U\, ;}
but it is well-defined globally on $S$, as the local definitions
coincide on intersections of different $S\cap U$'s by the
second point in Theorem 4.1.}
\medskip
\noindent
{\sc Remark 1}. \quad Note that $\res[\phi]\in
\Cal E^{p-1}(S)$ and is a closed form on $S$. Define the closed
set $\singsupp\left|_S\right.\phi$ on $S$ by saying that
$x\notin\singsupp\left|_S\right.\phi$ iff there is an
open neighborhood $U$ of $x$ in $M$ such that $\phi\left|_{U\setminus S}
\right.$ is bounded. Then from the above discussion we obtain:
\form{\supp\,\res[\phi]
\subset\singsupp\left|_S\phi\right.\, .}
\medskip
\noindent
{\sc Remark 2.}\quad When $k=0$, $M$ is a complex manifold
of complex dimension $n$ (Newlander-Nirenberg theorem) and $S$
is a complex hypersurface. In this case we recover exactly the
{\it forme-r\'esidu} of Leray [Lr]; which in turn was a generalization
of the classical residue in one complex variable.

\se{Properties of the residue form}
Let $\phi\in\Cal E^p(M\setminus S)$ be closed in $M\setminus S$
and have a first order pole along $S$; we investigate some
basic properties of $\res[\phi]$.
\prop{Let $\chi\in\Cal E^q(M)$ be a closed form on $M$. Then
$\phi\wedge\chi$ has a first order pole along $S$, and
\form{\res[\phi\wedge\chi]=\res[\phi]\wedge\chi\left|_S\right.\, .}}
\dimo If $\phi=\dfrac{ds}{s}\wedge\psi+\theta$, then
$\phi\wedge\chi=\dfrac{ds}{s}\wedge{\left(\psi\wedge\chi\right)}
+\left(\theta\wedge\chi\right)$, so the result is a consequence
of the uniqueness point of Theorem 4.1.
\medskip
Suppose that $S$ has a global smooth defining function $s$:
$$S\, = \, \left\{x\in M\left| s(x)=0\right.\right\}\, ,$$
with $ds\neq 0$ and $\bar\partial_Ms=0$ on $S$. Let
$\phi_1\in\Cal E^{p_1}(M\setminus S)$ and $\phi_2\in
\Cal E^{p_2}(M\setminus S)$ both be closed in $M\setminus S$,
and such that $s\phi_1$ and $s\phi_2$ are smooth across $S$.
\prop{Under the above assumptions:
\roster
\item"($i$)" $\phi_1\wedge\phi_2$ has a first order pole along $S$.
\item""
\item"($ii$)" $\phi_1\wedge\dfrac{ds}{s}$ and
$\phi_2\wedge\dfrac{ds}s$ also have first order poles along $S$.
\item""
\item"($iii$)" We have:
\form{\res\left[\phi_1\wedge\phi_2\right]=
\res\left[\phi_1\wedge\dfrac{ds}s\right]\wedge\res\left[\phi_2\right]
+\res\left[\phi_1\right]\wedge\res\left[\dfrac{ds}s\wedge\phi_2\right]\, .}
\endroster}
\dimo As in the previous discussion, we have
$$\phi_1=\dfrac{ds}s\wedge\psi_1+\theta_1
\qquad\text{and}\qquad
\phi_2=\dfrac{ds}s\wedge\psi_2+\theta_2$$
with $\psi_1,\, \psi_2,\, \theta_1, \, \theta_2$
smooth across $S$. Hence
\form{\matrix\format\l\\
\phi_1\wedge\phi_2\, = \,
\dfrac{ds}s\wedge\left\{
\psi_1\wedge\theta_2+(-1)^{p_1}\theta_1\wedge\psi_2\right\}
+\theta_1\wedge\theta_2\, ,\\
\phi_i\wedge\dfrac{ds}s\, = \, \theta_i\wedge\dfrac{ds}s\, ,\quad
i=1,2\, .
\endmatrix}
From these two formulas we obtain ($i$), ($ii$) and that
$\res\left[\dfrac{ds}s\wedge\phi_i\right]={\theta_i}\left|_S\right.$, and
(5.2) follows from (5.3).
\medskip
Next we derive a pointwise estimate for $\res[\phi]$: choose
a smooth Hermitian metric on the fibers of $\Bbb CT^*M$; it
induces a corresponding Hermitian metric on the fibers of
$\Lambda^p\Bbb CT^*M$. We define a smooth Hermitian metric
on the fibers of $\Bbb CT^*S$ by the natural identification
with the orthogonal complement of $\Bbb CT^*_SM$ in the
restriction to $S$ of $\Bbb CT^*M$; thereby inducing a
corresponding Hermitian metric on the fibers of
$\Lambda^{p-1}\Bbb CT^*S$. Let $\phi\in\Cal E^p(M\setminus S)$
be closed in $M\setminus S$ and have a first order pole
along $S$. Let $x\in S$ and $s$ be a smooth defining function
for $S$ near $x$, with $\bar\partial_Ms=0$ on $S$, and $s\phi$
be smooth across $S$ near $x$. Then we have
\prop{\form{
\left|\res[\phi](x)\right|_S\leq
\dfrac{\left|(s\phi)(x)\right|_M}{\left|ds(x)\right|_M}\, ,}
where $|\quad|_M$ and $|\quad|_S$ denote the Hermitian
lenghts in $\Lambda^p\Bbb CT^*_xM$, $\Lambda^1\Bbb CT^*_xM$ and
$\Lambda^{p-1}\Bbb CT^*_xS$, respectively.}
\dimo We note that for any form $\psi$ on $M$ we have
$\left|ds(x)\right|_M\cdot
\left|\left(\psi\left|_S\right)(x)\right.\right|_S\leq
\left|\left(ds\wedge\psi\right)(x)\right|_M$.
Since $s\phi=ds\wedge\psi+s\theta$, the result follows.
\medskip
We now study the behaviour of $\res[\phi]$ under the
pullback by a $CR$ mapping: Let $M^*$ be a smooth abstract
$CR$ manifold of type $(n^*,k^*)$, with the same hypotheses
as on $M$ (i.e. connected, paracompact, orientable). Let
$F:M^*@>>>M$ be a smooth $CR$ map such that $F(M^*)\cap S\neq
\emptyset$. We assume that {\sl $S$ is $CR$ transversal to $S$}.
By this we mean first of all that $F$ is transversal to $S$
in the sense of differential topology; i.e., that
\form{dF(x^*)\left(T_{x^*}M^*\right)+T_xS\, = \, T_xM\, ,\qquad
\forall x^*\in F^{-1}(S)\, ,\quad
x=F(x^*)\, ,}
and in addition that
\form{\cases
\dim dF(x^*)\left(H_{x^*}M^*\right)\quad\text{is constant for
$x^*\in F^{-1}(S)$}\\
\qquad\text{and}\\
dF(x^*)\left(H_{x^*}M^*\right)+H_xS\,=\, H_xM\, ,\qquad
\forall x^*\in F^{-1}(S)\, ,\quad
x=F(x^*)\, .
\endcases}
Let $S^*=F^{-1}(S)$. By (5.5) we have that $S^*$ is a smooth
closed submanifold of $M^*$ of real codimension two. From (5.6)
it follows that $S^*$ is a $CR$ submanifold of $M^*$, of type
$(n^*-1,k^*)$,
polar in $M^*$. Indeed let $x^*\in S^*$ and $x=F(x^*)$. Let
$s$ be a local defining function for $S$ near $x$, with
$\bar\partial_Ms=0$ along $S$. Then $s^*=s\circ F$ is a local
defining function for $S^*$ near $x^*$. Since $F$ is a $CR$ map,
we may replace $H$ in (5.6) by $T^{1,0}$ or $T^{0,1}$. We
first show that $\bar\partial_{M^*}s^*=0$ along $S^*$ near
$x^*$. Consider a point $y^*\in S^*$ near $x^*$; we have
$$\left\langle ds^*(y^*),T^{0,1}_{y^*}M^*\right\rangle
\, = \, \left\langle ds(y),dF(y^*)\left(T^{0,1}_{y^*}M^*\right)
\right\rangle=0\, ,\quad y=F(y^*)\, ,$$
because $dF(y^*)\left(T^{0,1}_{y^*}M^*\right)\subset
T^{0,1}_yM$. Finally we show that $S^*$ is polar in $M^*$.
From (5.6) we get
$$\matrix\format\l\\
\left\langle ds^*(y^*),T_{y^*}^{1,0}M^*\right\rangle=
\left\langle ds(y),dF(y^*)\left(T^{1,0}_{y^*}M^*\right)\right\rangle\\
 \quad \\
=\left\langle ds(y),dF(y^*)\left(T^{1,0}_{y^*}M^*\right)\right\rangle\\
 \quad\\
=\left\langle ds(y),T^{1,0}_yM\right\rangle\, \neq \, 0\, .
\endmatrix$$
This means that $ds^*(y^*)\in\Cal I_{y^*}(M^*)\setminus
\overline{\Cal I_{y^*}(M^*)}$, which forces $S^*$ to be
polar in $M^*$.
\par
Now consider a form $\phi\in\Cal E^p(M\setminus S)$, closed
in $M\setminus S$, and having a first order pole along $S$.
Then $F^*\phi\in\Cal E^p(M^*\setminus S^*)$, $F^*\phi$ is
closed in $M^*\setminus S^*$, and again has a first order
pole along $S^*$, since $F^*(s\phi)=s^*F^*\phi$.
\prop{In the situation described above, we have that
\form{\res\left[F^*\phi\right]\, = \, F^*\res[\phi]\, .}}
\dimo From $\phi=\dfrac{ds}s\wedge\psi+\theta$ it follows
that $F^*\phi=\dfrac{ds^*}{s^*}\wedge F^*\psi+F^*\theta$;
hence $\left(F^*\psi\right)\left|_{S^*}\right.=
F^*\left(\psi\left|_S\right.\right)$. [Here for simplicity
we use the same letter $F$ to denote the restriction map:
{$S^*\ni y^* @>>> F(y^*)\in S$}.]
\medskip
Finally consider a smooth closed submanifold $\Sigma_1$ of $M$
which is transversal to~$S$.
\prop{If $\phi\left|_{\Sigma_1\setminus S}\right.=0$,
then $\res[\phi]\left|_{\Sigma_1\cap S}\right.=0$.}
\dimo Locally we have $s\phi=ds\wedge\psi+s\theta$. Taking
pullbacks to $\Sigma_1$ we obtain
$0=(s\phi)\left|_{\Sigma_1}\right.=(ds)\left|_{\Sigma_1}\right.
\wedge\psi\left|_{\Sigma_1}\right.+(s\theta)\left|_{\Sigma_1}\right.$,
which along $\Sigma_1\cap S$ yields
$(ds)\left|_{\Sigma_1}\right.\wedge\psi\left|_{\Sigma_1}\right.=0$.
By the transverality assumption, $(ds)\left|_{\Sigma_1}\right.\neq 0$.
Therefore by Cartan's lemma, locally, $\psi\left|_{\Sigma_1}\right.
=(ds)\left|_{\Sigma_1}\right.\wedge\alpha+s\beta+\bar s\gamma$
for smooth forms $\alpha,\, \beta,\, \gamma$. Hence
$\left(\psi\left|_S\right.\right)\left|_{\Sigma_1}\right.
=\left(\psi\left|_{\Sigma_1}\right.\right)\left|_S\right.=0$
as claimed.

\se{The residue formula}
With $M$, $S$ as in the beginning of \S 3, we introduce also
\form{\Sigma\, = \Sigma_1\, \cup \, \Sigma_2\, \cup \,
\cdots \, \cup \Sigma_\ell\, ,}
where each $\Sigma_j$ is a smooth closed submanifold of $M$.
We assume that the $S$, $\Sigma_1$, $\Sigma_2$, $\hdots$,
$\Sigma_\ell$ are {\it in general position}.\par
This implies, in particular, that the intersection of any
subset of the $\Sigma_1, \, \hdots , \, \Sigma_\ell$
is transversal to $S$. This allows to constuct a tubular
neighborhood $V$ of $S$ in $M$ adapted to $\Sigma$: recall
that the normal bundle $N_SM$ of $S$ in $M$ is the
quotient $ TM\left|_S\right.\left/_{\dsize{TS}}\right.$.
The tubular neighborhood $V$ is the datum of an open neighborhood
$V$ of $S$ in $M$, together with a smooth diffeomorphism
$$\CD
\nu: N_SM @>{\sim}>> V
\endCD$$
which is the identity on $S$, identified with the zero section
of $N_SM$. Denote by $\pi:N_SM\longrightarrow S$ the projection.
The fact that the tubular neighborhood is {\sl adapted} to $\Sigma$
means, in particular, that $\nu\left(\pi^{-1}(\Sigma\cap S)\right)
=\Sigma\cap V$. \par
This inables us to construct a smooth strict deformation retract
$$\mu:[0,1] \times V \longrightarrow V$$
of the pair $(V,\Sigma\cap V)$ onto the pair $(S,\Sigma\cap S)$;
namely,
\form{\mu(t,x)\, = \, \mu_t(x)\, = \, \nu\left(t\,\nu^{-1}(x)\right)\, .}
We fix a smooth Riemannian metric on the vector bundle
$N_SM$, and denote the lenght of a vector $v$ in the fiber at $x\in S$
by $|v|_x$.\par
Consider a relative $(p-1)$-cycle $\gamma$ in $(S,\Sigma\cap S)
\equiv (S,\Sigma)$. Here we use smooth singular chains with compact
support, and take $\Bbb Z$-coefficients. It is standard to identify
$\gamma$ with a smooth map
\form{\tilde\gamma:(P_\gamma,\partial P_\gamma)
\longrightarrow (S,\Sigma)\, ,}
where $P_\gamma$ is a finite polyhedron, of dimension $(p-1)$,
embedded in some Euclidean space $\Bbb R^N$. Let
$\hat\gamma^*\left(N_SM\right)$ denote the pulled-back bundle
over $P_\gamma$. Using the metric on the fibers, we introduce the
closed disk bundle $D_\gamma(t)$ and the circle bundle
$C_\gamma(t)$ of radius $t>0$ contained in $\hat\gamma^*\left(
N_SM\right)$. These can be regarded as $(p+1)$ and $p$ dimensional
polyhedra, respectively, embedded in some $\Bbb R^{N'}$.
The map $\hat\gamma$ lifts to a smooth vector bundle
homomorphism
$$\hat\gamma_{N_SM}:\hat\gamma^*\left(N_SM\right)
\longrightarrow N_SM\, .$$
This enables us to define
\form{\cases
\widehat{D_t\gamma}:D_\gamma(t)\longrightarrow V \subset M\, ,
\\ \qquad
\\
\qquad \text{and}\\ \quad
\\
\widehat{\delta_t\gamma}:C_\gamma(t)\longrightarrow V\setminus S
\subset M\setminus S\, ,
\endcases}
by $\left(\widehat{D_t\gamma}\right)(y)=\nu\left(\hat\gamma_{N_SM}(y)\right)$
for $y\in D_\gamma(t)$, and
$\left(\widehat{\delta_t\gamma}\right)(y)
=\nu\left(\hat\gamma_{N_SM}(y)\right)$
for $y\in C_\gamma(t)$.
Since
\form{\cases
\partial C_\gamma(t)\, = \,
{C_\gamma(T)}\left|_{\partial P_\gamma}\right.\, ,\\
\partial D_\gamma(t)\, = \, C_\gamma(t)\, \cup\,
{D_\gamma(t)}\left|_{\partial P_\gamma}\right.
\endcases}
and our tubular neighborhood $V$ was adapted to $\Sigma$,
$\widetilde{\delta_t\gamma}$ defines a smooth singlular
relative $p$-cycle $\delta_t\gamma$ in $(M\setminus S,\Sigma)$,
and $\widetilde{D_t\gamma}$ defines a smooth singular relative
${(p+1)}$-chain $D_t\gamma$ in $(M,\Sigma)$, such that
$\delta_t\gamma=\partial'\left(D_t\gamma\right)$, where
$\partial'$ denotes the relative boundary in $(M,\Sigma)$,
for each $0<t\leq 1$. Note that each $\delta_t\gamma$
is homologous to $\delta_1\gamma$ in
$(M\setminus S,\Sigma)$ by the construction. \par
We now turn to the theorem establishing the
residue formula. We fix the following notation:
$\Cal E^p(M,\Sigma)$, $\Cal E^p(M\setminus S,\Sigma)$ and
$\Cal E^{p-1}(S,\Sigma)$ will denote the subspaces of
$\Cal E^p(M)$, $\Cal E^p(M\setminus S)$ and $\Cal E^{p-1}(S)$
consisting of those smooth forms whose pullbacks
to $\Sigma$, $\Sigma\cap\left( M\setminus S\right)$ and
$\Sigma\cap S$ are zero, respectively. \par
Let $\gamma$ be a smooth compactly supported
singular relative $(p-1)$-cycle in $(S,\Sigma)$.
Likewise let $\Gamma$ be a smooth compactly supported
singular relative $p$-cycle in $(M\setminus~S,\Sigma)$.
If $\Gamma$ is homologous to $\delta_1\gamma$ in
$(M\setminus S,\Sigma)$, we shall say that
{\it $\Gamma$ is cobordant to $\gamma$ in $(M\setminus S,\Sigma)$}.
\thm{Let $\phi\in\Cal E^p(M\setminus S,\Sigma)$ be closed
in $M\setminus S$ and have a first order pole along
$S$. Let $\gamma$ be a smooth compactly supported singlular
relative $(p-1)$-cycle in $(S,\Sigma)$. Then for every
$\Gamma$ wich is "cobordant to $\gamma$ in $(M\setminus S,\Sigma)$",
we have
\form{\dsize\int_{\Gamma}{\phi}\, = \,
2\pi\ssize{\sqrt{-1}}\, \dsize\int_{\gamma}{\res[\phi]}\, .}}
\dimo
By Stokes' theorem we obtain
\form{\dsize\int_{\Gamma}{\phi}\, = \,
\dsize\int_{\delta_1\gamma}{\phi}\, = \,
\dsize\int_{\delta_t\gamma}{\phi}\qquad
\text{for}\quad 0<t\leq 1\, ,}
as $\Gamma\sim\delta_1\gamma\sim\delta_t\gamma$ in
$(M\setminus S,\Sigma)$.
Hence it will suffice to show that
\form{\lim_{t \searrow 0}{\dsize\int_{\delta_t\gamma}{\phi}}
\, = \, 2\pi\ssize{\sqrt{-1}}\,\dsize\int_{\gamma}{\res[\phi]}\, ,}
which is a limit already known to exist, by (6.7).
\par
Along $S$ we introduce a locally finite partition of unity
$\sum_{\alpha}{\chi_\alpha}\equiv 1$, subordinate to a covering
$\left\{U_\alpha\right\}$ of $S$ by open neighborhoods $U_\alpha$
in which $S\cap U_\alpha$ has a locally defining funciton
$s_\alpha$ satisfying $\bar\partial_Ms_\alpha=0$ along $S$,
and such that $s_\alpha\phi\in \Cal E^p(U_\alpha)$. Then
we have the decomposition
\form{\dsize\int_{\delta_t\gamma}{\phi}
\, = \,
\dsize\int_{\delta_t\gamma}{\left(\sum_\alpha{\chi_\alpha\phi}\right)}
\, = \,
\dsize\sum_\alpha\dsize\int_{\delta_t\gamma}{\chi_\alpha\phi}\, ,}
which can be regarded as a finite sum, uniformly for $0<t\leq 1$.
Hence it will suffice to show that for each $\alpha$,
\form{\lim_{t\searrow 0}{\dsize\int_{\delta_t\gamma}{\chi_\alpha\phi}}
\, = \,
2\pi\ssize{\sqrt{-1}}\,\dsize\int_{\gamma}{\chi_\alpha\res[\phi]}\, .}
From Theorem 4.1 we have
$$\phi=\dfrac{ds_\alpha}{s_\alpha}\wedge\psi+\theta\quad\text{in}\quad
U_\alpha\, ,$$
with $\psi$ ane $\theta$ smooth across $S$. By Stokes' theorem
$$\left| \dsize\int_{\delta_t\gamma}{\chi_\alpha\theta}\right|
\, = \,
\left| \dsize\int_{D_\gamma(t)}{d(\chi_\alpha\theta)}\right|\,
\leq\, \text{const}\cdot t^2 @>>> 0\, ,$$
as $t\searrow 0$, since $\theta$ is smooth across $S$.
Using the retraction $\mu_0:V @>>> S$ from (6.2), we may write
\form{\psi-\mu_0^*\left(\psi\left|_S\right.\right)\, = \,
s_\alpha\psi_1+\bar s_\alpha\psi_2+ds_\alpha\wedge\psi_3+
d\bar s_\alpha\wedge\psi_4}
in $U_\alpha$, with smooth forms $\psi_i$, because the
left hand side in (6.11) pulls back to zero along $S$. Note that
$\left[\mu_0^*\left(\psi\left|_S
\right.\right)\right]\left|_{\Sigma\cap V}\right.=0$ by
Proposition 5.5.
Consider
\form{\dfrac{ds_\alpha}{s_\alpha}\wedge\left[
\psi-\mu^*_0\left(\psi\left|_S\right.\right)\right]\, = \,
ds_\alpha\wedge\psi_1+\dfrac{\bar s_\alpha}{s_\alpha}d\bar s_\alpha
\wedge\psi_2+\dfrac{ds_\alpha\wedge d\bar s_\alpha}{s_\alpha}
\wedge\psi_4\, .}
Note that each term on the right in (6.12) is uniformly bounded
on the intersection of any compact subset of $U_\alpha$ with
$U_\alpha\setminus S$. Hence
$$\left| \dsize\int_{\delta_t\gamma}{
\dfrac{ds_\alpha}{s_\alpha}\wedge\left[\psi-
\mu^*_0\left(\psi\left|_S\right.\right)\right]\,\chi_\alpha}
\right| \, \leq \, \text{const}\cdot t @>>> 0\, ,$$
as $t\searrow 0$. It follows that
$$\lim_{t\searrow 0}{\dsize\int_{\delta_t\gamma}{\chi_\alpha\phi}}
\, = \,
\lim_{t\searrow 0}{\dsize\int{\chi_\alpha\,\dfrac{ds_\alpha}{s_\alpha}
\wedge \mu_0^*\left(\psi\left|_S\right.\right)}}
\, = \,
2\pi\ssize{\sqrt{-1}}\,\dsize\int_\gamma{\chi_\alpha\,\psi\left|_S\right.}
\, ,$$
and the theorem is proved.
\medskip \noindent
{\sc Remark.}\quad As we are using smooth singular relative
homology with $\Bbb Z$-coefficients, it may happen that the
relative cycle $\gamma$ has finite order, say $a\cdot \gamma
\sim 0$ in $\ho_{p-1}(S,\Sigma;\Bbb Z)$ for some integer $a\neq 0$.
In that case both right and left hand sides in the
residue formula (6.6) are zero. Indeed the right-hand
side is zero by Stokes' theorem, since $\res[\phi]$ is
closed on $S$ and pullbacks to zero on $S\cap \Sigma$;
likewise the left-hand side is zero because $a\Gamma \sim
a\delta_1\gamma\sim 0$ in $\ho_p(M\setminus S,\Sigma;\Bbb Z)$.

\se{The Homological Residue}
In what follows we shall denote the standard relative homology
groups, with compact support and integer coefficients, of $(M,\Sigma)$,
$(S,\Sigma)$ and $(M\setminus S,\Sigma)$ by
$\ho_p(M,\Sigma;\Bbb Z)$, $\ho_p(S,\Sigma;\Bbb Z)$ and
$\ho_p(M\setminus S,\Sigma;\Bbb Z)$, respectively.
As $M$, $S$ are smooth we may compute these, as is well-known,
using {\it smooth} singular relative $p$-chains.\par
As a consequence of the discussion in \S 6, we actually obtain
a {\it coboundary} homomorphism
\form{\delta:\ho_{p-1}(S,\Sigma;\Bbb Z)\longrightarrow
\ho_p(M\setminus S,\Sigma;\Bbb Z)\, ,}
defined by $\delta\left([\gamma]\right)=\left[\delta_1\gamma\right]$;
indeed the construction we have given of $\delta_1\gamma$ commutes
with the boundary maps. This passes to a homomorphism
$$\bar{\delta}:\bar{\ho}_{p-1}(S,\Sigma;\Bbb Z)
\longrightarrow \bar{\ho}_p(M\setminus S,\Sigma;\Bbb Z)$$
of the {\it weak} homology groups. Here we use a bar to indicate
the weak homology groups defined by
$\bar{\ho}_p=\ho_p\left/_{\dsize{\roman{Tor}_p}}\right.$.
Therefore by Stokes' theorem the residue formula (6.6) can
be reformulated more generally in terms of weak homology:
\thm{Let $\phi\in\Cal E^p(M\setminus S,\Sigma)$ be closed
in $M\setminus S$ and have a first order pole along $S$.
Let $\bar{h}\in\bar{\ho}_{p-1}(S,\Sigma;\Bbb Z)$
be a coset of relative homology classes. Then
\form{\dsize\int_{\bar{\delta}\bar{h}}{\phi}\, = \,
2\pi\ssize{\sqrt{-1}}\, \dsize\int_{\bar{h}}{\res[\phi]}\, ,}
in which $\bar\delta\bar{h}$ is a coset of relative homology
classes in $\bar{\ho}_p(M\setminus S,\Sigma;\Bbb Z)$.}
\medskip
We now return to the "coboundary" homomorphism $\delta$ in (7.1);
it can be inserted into a long exact sequence as follows: Set
$$\nu\left(\left\{\xi\in N_SM\,\left|\,|\xi|_{\pi(\xi)}\leq 1
\right.\right\}
\right)\, = \, V_1\, \subset V$$
and consider the standard long exact sequence of a relative pair
\form{\matrix\format\l\\
\cdots\rightarrow
\ho_{p+1}\left(\overline{M\setminus V_1},\Sigma;\Bbb Z\right)
\overset{i}\to\longrightarrow
\ho_{p+1}(M,\Sigma;\Bbb Z) \\
\quad\\
\qquad\qquad
\overset{p}\to\longrightarrow
\ho_{p+1}\left(M,\overline{M\setminus V_1}\cup\Sigma;\Bbb Z\right)
\overset\partial\to\longrightarrow
\ho_p(\overline{M\setminus V_1},\Sigma;\Bbb Z)
\rightarrow\cdots
\endmatrix}
First we observe that
$$\ho_p\left(\overline{M\setminus V_1},\Sigma;\Bbb Z\right)
\simeq \ho_p(M\setminus S,\Sigma;\Bbb Z)$$
because
$\left(\overline{M\setminus V_1},\Sigma\right)$ is
a relative deformation retract of $(M\setminus S,\Sigma)$.
Next we obtain
$$\ho_{p+1}\left(M,\overline{M\setminus V_1}\cup\Sigma;
\Bbb Z\right)\simeq \ho_{p+1}(V_1,\partial V_1\cup\Sigma;\Bbb Z)$$
by the excision of $M\setminus V_1$.	But
by the relative Thom isomorphism we have
$$\ho_{p+1}(V_1,\partial V_1\cup\Sigma;\Bbb Z)\simeq
\ho_{p-1}(S,\Sigma;\Bbb Z)\, .$$
Thus we arrive at (cf. Leray [Lr]):
\prop{There is a long exact sequence
\form{\matrix\format\l\\
\cdots @>>> \ho_{p+1}(M\setminus S,\Sigma;\Bbb Z)
\longrightarrow
\ho_{p+1}(M,\Sigma;\Bbb Z) \\
\quad\\
\qquad\qquad
\overset\tau\to\longrightarrow
\ho_{p-1}(S,\Sigma;\Bbb Z)
\overset\delta\to\longrightarrow
\ho_p(M\setminus S,\Sigma;\Bbb Z) \rightarrow\cdots
\endmatrix}}
\medskip
\noindent
{\sc Remark.}\quad The connecting homomorphism $\partial$
in (7.3) becomes our $\delta$ in (7.4) because the
inverse of the Thom isomorphism is the map
$[\gamma]@>>>\left[D_1\gamma\right]$ described in \S 6,
which becomes $[\gamma]@>>>\left[\delta_1\gamma\right]$
when composed with the boundary map $\partial$.

\se{The Cohomological Residue}
In what follows we assume $M$, $S$ as in the beginning of \S 3, and
$\Sigma$ as in the beginning of \S 6. Since $d$ commutes with the
pullback to $\Sigma$, we have the complex:
\form{\CD
0\rightarrow \Cal E^0(M,\Sigma) @>d>> \Cal E^1(M,\Sigma)@>d>>
\Cal E^2(M,\Sigma) \rightarrow\cdots
\endCD}
This is the de Rham complex on smooth (complex valued) forms $\phi$
on $M$ with $\phi\left|_{\Sigma}\right.=0$. (We say that $\phi$
has {\it zero Cauchy data} for $d$ along $\Sigma$). We denote the
$p$-th cohomology group of (8.1) by $\ho^p(M,\Sigma)$. Instead of the
pair $(M,\Sigma)$ we may use the pair
$(S,\Sigma)\equiv (S,S\cap\Sigma)$, or the pair $(M\setminus S,
\Sigma)\equiv \left(M\setminus S,(M\setminus S)\cap\Sigma\right)$
in (8.1). So we have also $\ho^p(S,\Sigma)$ and $\ho^p(M\setminus S,\Sigma)$.
Since these cohomology groups are not entirely standard, we
relate them to the standard singular relative cohomology
groups $\ho^p_{\text{sing}}$ with complex coefficients:
\prop{For all $p$
\roster
\item"($a$)" $\ho^p(M,\Sigma)\simeq \ho^p_{\text{sing}}(M,\Sigma;\Bbb C)$
\item"($b$)" $\ho^p(S,\Sigma)\simeq \ho^p_{\text{sing}}(S,\Sigma;\Bbb C)$
\item"($c$)" $\ho^p(M\setminus S,\Sigma)\simeq
\ho^p_{\text{sing}}(M\setminus S,\Sigma;\Bbb C)$
\endroster}
\dimo We do the proof for ($a$), as ($b$) and ($c$) are the same.
[Indeed, we can replace $M$ by any smooth submanifold $\Omega$ of
$M$, provided that $\Omega$, $\Sigma_1$, $\Sigma_2$, $\hdots$,
$\Sigma_\ell$ are in general position.] Note that when $\Sigma=
\emptyset$, these are just the isomorphisms given by de Rham's
theorem.\par
First suppose $\Sigma=\Sigma_1$. For each $p$ we have the exact
sequence:
\form{\CD
0 \rightarrow
\Cal E^p(M,\Sigma) @>>> \Cal E^p(M) @>>>\Cal E^p(\Sigma)
\rightarrow 0\, .
\endCD}
Since the maps in (8.2) commute with $d$, we obtain the long
exact sequence:
\form{\CD
\cdots\rightarrow
\ho^{p-1}(\Sigma) @>>> \ho^p(M,\Sigma) @>>> \ho^p(\Sigma)
\rightarrow\cdots
\endCD}
We consider the corresponding sequence in singular relative
cohomology:
\form{\CD
\cdots\rightarrow
\ho^{p-1}_{\roman{sing}}(\Sigma;\Bbb C)
@>>> \ho^p_{\roman{sing}}(M,\Sigma;\Bbb C) @>>>
\ho^p_{\roman{sing}}(\Sigma;\Bbb C)
\rightarrow\cdots
\endCD}
By de Rham's theorem we have that
$$\ho^p(M)\simeq \ho^p_{\roman{sing}}(M;\Bbb C)\quad\text{and}\quad
\ho^p(\Sigma)\simeq \ho^p_{\roman{sing}}(\Sigma;\Bbb C)\, .$$
As these isomorphisms are compatible with the natural
homomorphisms between (8.3) and (8.4), we obtain that also
$\ho^p(M,\Sigma)\simeq \ho^p_{\roman{sing}}(M,\Sigma;\Bbb C)$
upon application of the $5$-lemma. \par
We proceed by induction on $\ell$: Assume that the result
is known for $\Sigma'=\Sigma_1\cup\Sigma_2\cup\cdots\cup
\Sigma_{\ell-1}$, and consider $\Sigma=\Sigma'\cup\Sigma_{\ell}$.
By our transversality assumption, there is for each $p$ a
short exact sequence
\form{
0\rightarrow
\Cal E^p(M,\Sigma)
\longrightarrow
\Cal E^p(M,\Sigma')\oplus\Cal E^p(M,\Sigma_\ell)
\longrightarrow
\Cal E^p(M,\Sigma'\cap\Sigma_\ell)
\rightarrow 0\, .
}
As the maps in (8.5) commute with $d$, we obtain the long exact
sequence:
\form{\matrix\format\l&\l&\r&\r&\l&\l\\
&&\cdots&@>>>&
\ho^{p-1}(M,\Sigma'\cap\Sigma_\ell)
&\rightarrow\\
\rightarrow
\ho^p(M,\Sigma)
&@>>>&
\ho^p(M,\Sigma')\oplus \ho^p(M,\Sigma_\ell)
&@>>>& \ho^p(M,\Sigma'\cap\Sigma_\ell)&\rightarrow\cdots
\endmatrix}
Since $\Sigma'\cap\Sigma_\ell=(\Sigma_1\cap\Sigma_\ell)\cup
\cdots\cup(\Sigma_{\ell-1}\cap\Sigma_\ell)$ is a union of
$(\ell-1)$ smooth submanifolds in general position, we can
apply our inductive hypothesis to conclude that
$$\ho^p(M,\Sigma'\cap\Sigma_\ell)\simeq \ho^p_{\roman{sing}}(M,
\Sigma'\cap\Sigma_\ell;\Bbb C)\, .$$
Hence we may compare (8.6) with
the corresponding Mayer-Vietoris sequence in singular
relative cohomology, and apply the $5$-lemma once again to
obtain $\ho^p(M,\Sigma)\simeq \ho^p_{\roman{sing}}(M,\Sigma;\Bbb C)$.
This completes the proof.
\medskip
\prop{There is a long exact sequence:
\form{
\cdots\rightarrow
\ho^p(M,\Sigma)@>>>
\ho^p(M\setminus S,\Sigma)
@>{\delta^*}>>
\ho^{p-1}(S,\Sigma)
@>{\tau^*}>>
\ho^{p+1}(M,\Sigma)
@>>>
\cdots\, .}}
\dimo By the universal coefficient theorem we have that
$$\ho^p_{\roman{sing}}(M,\Sigma;\Bbb C)\simeq
\Hom_{\Bbb Z}\left(\ho_p(M,\Sigma;\Bbb Z),\Bbb C\right)\, ,$$
and likewise for the other pairs $(S,\Sigma)$ and
$(M\setminus S,\Sigma)$. Applying the
$\Hom_{\Bbb Z}(\,\cdot\, ,\Bbb C)$ functor to (7.4) we obtain
(8.7) because $\Bbb C$ is an injective $\Bbb Z$-module.
\medskip
\noindent
{\sc Definition}.{\sl \quad For $h^*\in \ho^p(M\setminus S,\Sigma)$,
\form{\Res h^*\, = \, \dfrac{1}{2\pi\ssize{\sqrt{-1}}}\,
\delta^* h^*\in \ho^{p-1}(S,\Sigma)}
will be called the {\bf residue class} of $h^*$.}
\medskip
\prop{Let $\phi\in\Cal E^p(M\setminus S,\Sigma)$ be closed in
$M\setminus S$ and have a simple pole along $S$. Then
\form{\Res[\phi]\, = \, \left[\res[\phi]\right]\, .}}
\dimo $$\matrix\format\r&&\l\\
\Res[\phi]\left([\gamma]\right)&=&\dfrac{1}{2\pi\ssize{\sqrt{-1}}}
[\phi]\left([\delta\phi]\right)\\
\quad\\
\quad\\
=\dfrac{1}{2\pi\ssize{\sqrt{-1}}}\dsize\int_{\delta\gamma}{\phi}
&=&\dsize\int_{\gamma}{\res[\phi]}\\
\quad\\
\quad\\
&=&\left[\res[\phi]\right]\left([\gamma]\right)\, ,
\\
\quad
\endmatrix$$
for every smooth singlular relative $(p-1)$-cycle in $(S,\Sigma)$.
So $\Res[\phi]$ and $\left[\res[\phi]\right]$ represent the same
element of
$\Hom_{\Bbb Z}\left(\ho_{p-1}(S,\Sigma;\Bbb Z),\Bbb C\right)$,
and the result follows from the discussion above.
\thm{Each cohomology class $h^*\in\ho^p(M\setminus S,\Sigma)$
has a representative $\phi\in\Cal E^p(M\setminus S,\Sigma)$,
closed in $M\setminus S$, and having a simple pole along $S$.}
\dimo
First we observe that $HM$ can be viewed as a {\sl complex}
smooth vector bundle of rank $n$, with complex structure on
the fibers given by $J$. Since $S$ is polar, $HS$ is a
smooth complex subbundle of $HM$ along $S$, having rank $(n-1)$.
Moreover the projection $TM\left|_S\right. @>>> N_SM$
gives by restriction a surjective vector bundle morphism
$HM\left|_S\right. @>>> N_SM$, which in turn yields an
isomorphism
\form{\CD
{HM\left|_S\right.}\left/_{\dsize{HS}}\right.
@>{\sim}>> N_SM \, .\endCD}
Thus $N_SM$ has the natural structure of a {\sl complex
line bundle} over $S$. We take a local trivialization
$\left\{(U_\alpha,s_\alpha)\right\}$ of $N_SM$ and use
$s_\alpha$ as the local defining functions for $S$ in
$\mu_0^{-1}(U_\alpha)$, with $\mu_0$ given in (6.2):
we denote by $s_\alpha(x)$ the value of the $s_\alpha$-coordinate
corresponding to $\nu^{-1}(x)$, for $x\in\mu_0^{-1}(U_\alpha)$.
Here $\nu:N_SM @>{\sim}>> V$ is the diffeomorphism from \S 6,
which represents $S$ as the zero section in $N_SM$. We observe that
$\bar\partial_Ms_\alpha=0$ along $S$ because the complex structure
on the fibers of $N_SM$ agrees with the partial complex
structure of $M$ along $S$. Moreover for
$x\in\mu_0^{-1}(U_\alpha)\cap\mu_0^{-1}(U_\beta)$, we have that
\form{s_\alpha(x)=g_{\alpha\beta}(\mu_0(x))s_\beta(x)\, ,}
where the $g_{\alpha\beta}$ are the transition functions,
$g_{\alpha\beta}\neq 0$, corresponding to the local trivialization
of the bundle. From (8.11) we derive that
\form{\dfrac{ds_\alpha}{s_\alpha}-\dfrac{ds_\beta}{s_\beta}
\, = \, \dfrac{dg_{\alpha\beta}}{g_{\alpha\beta}}\, .}
Since the right-hand side of (8.12) is smooth across $S$, it
follows that the left-hand side has a smooth extension
across $S$.\par
Now consider a class $h^*\in\ho^p(M\setminus S,\Sigma)$. Let
$\psi_0\in\Cal E^{p-1}(S,\Sigma)$ with $d\psi_0=0$ be
any representative of $\delta^* h^*$. We choose a smooth partition
of unity $\left\{\chi_\alpha\right\}$ of a neighborhood of $S$
in $V$, subordinated to the covering
$\left\{\mu_0^{-1}(U_\alpha)\right\}$, and with
$\supp \chi_\alpha\subset\nu\left(
\left\{\xi\in N_SM\,\left| \, |\xi|_{\pi(\xi)}\leq \frac{3}{4}
\right.\right\}\right)=V_{\frac{3}{4}}\subset V$.
Let
$$\omega\, = \, \dsize\sum_\alpha{\chi_\alpha\dfrac{ds_\alpha}{s_\alpha}
\wedge \mu_0^*(\psi_0)}\, .$$
Then $\omega\in \Cal E^p(M\setminus S,\Sigma)$ and $\omega=0$ outside
of $V_{\frac{3}{4}}$. It follows from (8.12) that $d\omega$ has a
smooth continuation across $S$, and defines a closed form in
$\Cal E^{p+1}(M,\Sigma)$: Indeed in a sufficiently small open
neighborhood of $S\cap\mu_0^{-1}(U_\beta)=U_\beta$ in $M$,
\form{\matrix\format\r\,&&\,\l\\
\omega&=&\dfrac{ds_\beta}{s_\beta}\wedge\mu_0^*(\psi_0)
+\dsize\sum_\alpha{\chi_\alpha\left(
\dfrac{ds_\beta}{s_\beta}
-\dfrac{ds_\alpha}{s_\alpha}
\right)\wedge\mu_0^*(\psi_0)}\\
&&\quad\\
&=&
\dfrac{ds_\beta}{s_\beta}
\wedge\mu_0^*(\psi_0)+
\dsize\sum_\alpha{\chi_\alpha
\dfrac{dg_{\alpha\beta}}{g_{\alpha\beta}}
\wedge\mu^*_0(\psi_0)}\, ;\endmatrix}
hence
\form{d\omega=\dsize\sum_\alpha{
d\chi_\alpha\wedge
\dfrac{dg_{\alpha\beta}}{g_{\alpha\beta}}\wedge
\mu^*_0(\psi_0)
}\, ,}
which is smooth across $S$. From formula (8.13) we see
that $\omega$ has a first order pole along $S$.\par
Next we show that $\tau^*\left([\psi_0]\right)=
-\dfrac{1}{2\pi\ssize{\sqrt{-1}}}\left[d\omega\right]$
in $\ho^{p+1}(M,\Sigma)$: consider a homology class
$\Omega\in\ho_{p+1}(M,\Sigma;\Bbb Z)$, and let $\gamma$
be a smooth singular relative $(p-1)$ cycle in $(S,\Sigma)$
representing the class
$\tau(\Omega)=[\gamma]\in\ho_{p-1}(S,\Sigma;\Bbb Z)$. Then
$\Omega-\left[D_1\gamma\right]$ is the image of a homology class
in $\ho_{p+1}(M\setminus S,\Sigma;\Bbb Z)$. By excision we can
find a smooth singular relative $(p+1)$-cycle $\Omega'$
in $\left(\overline{M\setminus V_1},\Sigma\right)$
such that $[\Omega']=\Omega-\left[D_1\gamma\right]$.
Then applying Stokes' theorem we obtain
$$\matrix\format\r\,&&\,\l\\
\left[d\omega\right](\Omega)&=&
\dsize\int_{\Omega}{d\omega}\, = \,
\dsize\int_{\Omega'}{d\omega}+\dsize\int_{D_1\gamma}{d\omega}\\
&&\\
&=&\dsize\int_{\delta_1\gamma}{\omega}-
\lim_{\epsilon\searrow 0}{
\dsize\int_{\delta_\epsilon\gamma}{\omega}}\\
&&\\
&=&-2\pi{\ssize{\sqrt{-1}}}\, \dsize\int_{\gamma}{\psi_0}\\
&&\\
&=& -2\pi{\ssize{\sqrt{-1}}}\, [\psi_0]\left(\tau(\Omega)\right)\\
&&\\
&=& -2\pi{\ssize{\sqrt{-1}}}\, \tau^*\left([\psi_0]\right)(\Omega)
&&\\
\, .
\endmatrix$$
We have that $\tau^*\left([\psi_0]\right)=\tau^*\delta^* h^*=0$;
hence $[d\omega]\sim 0$ in $\ho^{p+1}(M,\Sigma)$. Therefore
$d\omega=d\eta$ in $M\setminus S$, for some
$\eta\in\Cal E^p(M,\Sigma)$. Hence
$\omega-\eta$ is closed in $M\setminus S$, has a simple pole
along $S$, and its residue form is $\psi_0$. By
Proposition 8.3 we have that
$\delta^*\left(h^*-2\pi{\ssize{\sqrt{-1}}}[\omega-\eta]\right)
=0$, and therefore by the exact sequence (8.7) there exists
a closed form $\lambda\in\Cal E^p(M,\Sigma)$ such that
$h^*-2\pi{\ssize{\sqrt{-1}}}[\omega-\eta]=
\left[\lambda\left|_{M\setminus S}\right.\right]$.
Let $\phi=2\pi{\ssize{\sqrt{-1}}}(\omega-\eta)+\lambda$.
This is a closed form in $\Cal E^p(M\setminus S,\Sigma)$
having a simple pole along $S$, with
$[\phi]=h^*$ in $\ho^p(M\setminus S,\Sigma)$, completing
the proof of the theorem.
\medskip
We now come to the cohomological version of the residue
formula.
\thm{Let $h^*\in\ho^p(M\setminus S,\Sigma)$ and
$\bar h\in \bar{\ho}_{p-1}(S,\Sigma;\Bbb Z)$. Then
\form{\dsize\int_{\bar\delta\bar h}{h^*}
\, = \, 2\pi\ssize{\sqrt{-1}}\,\dsize\int_{\bar h}{\Res h^*}\, .}}
\dimo Note that both the left and the right-hands in (8.15)
are well defined by Stokes' theorem. The result is then a
consequence of Theorem 8.4 and Proposition 8.3.
\medskip
\se{Properties of the Residue Class and Global
defining Functions}
Next we summarize some important properties of the residue class
map which follow from the analogous properties of the
residue form, in view of Proposition 8.3 and Theorem 8.4.
\prop{If $h^*\in\ho^p(M\setminus S,\Sigma)$ and
$k^*\in\ho^q(M)$, then
\form{\Res\left( h^*\smallsmile k^*\left|_{M\setminus S}\right.\right)
\, = \, \left(\Res h^*\right)\smallsmile k^*\left|_{S}\right.}}
Here the cup $\smallsmile$ denotes the product operation on
cohomology induced by the wedge product of forms.
\prop{Assume that $S$ has a global smooth defining function
$s$ with $\bar\partial_Ms=0$ along $S$. If $h^*\in
\ho^p(M\setminus S,\Sigma)$ and $g^*\in\ho^q(M\setminus S,\Sigma)$
then
\form{\Res\left( h^*\smallsmile g^*\right) \, = \,
\Res\left(h^*\smallsmile\left[\dfrac{ds}{s}\right]\right)
\smallsmile\Res g^* \, + \,
\Res h^*\smallsmile \Res\left(\left[\dfrac{ds}{s}\right]
\smallsmile g^*\right)\, .}}
\medskip
Note that we may apply proposition 5.2. at this point, because
the construction in our proof of Theorem 8.4. allows us to find
representatives in each of the classes $h^*$ and $g^*$ having
a simple pole with respect to the same defining function $s$.
\smallskip
Next we consider the situation of a smooth $CR$ map
$F:M^*@>>>M$ as in \S 5. We assume that $F$ is $CR$ transversal
to $S$, and moreover that $F$ is transversal to each~$\Sigma_j$.
\prop{Under the above assumptions, we have
\form{\Res\left( F^*h^*\right)\, = \, F^*\left(\Res h^*\right)}
for any $h^*\in\ho^p(M\setminus S,\Sigma)$.}
\medskip
\noindent
{\sc Remark}.\quad In particular if $\Cal V$ is an open neighborhood of
$S$ in $M$, and $F$ is the inclusion map, we obtain
\form{\Res h^*\, = \, \Res\left( h^*\left|_{\Cal V}\right.\right)\, .}
Recall that the smooth complex line bundle $N_SM$ is trivial if
its Chern class in $\ho^2(S,\Bbb Z)$ is zero.
\prop{If the normal bundle $N_SM$ has zero Chern class, then
\roster
\item"($a$)"  $S$ has a global defining function $s$, with
$\bar\partial_M s=0$ along $S$, defined in a neighborhood $\Cal V$
of $S$ in $M$.
\item"($b$)" The formula analogous to (9.2) holds, in which $M$
is replaced by $\Cal V$.
\endroster}
\medskip
Finally we remark that, in general, by Proposition 9.1 we have that
$\Res$ is a homomorphism of $\ho^*(M)$-algebras from
$\ho^*(M\setminus S,\Sigma)$ to $\ho^*(S,\Sigma)$.
\medskip

\se{Higher Order Poles}
We say that $\left\{(U_\alpha,s_\alpha)\right\}$ form a {\it consistent}
system of local defining functions for $S$ iff: the $\{U_\alpha\}$
is a locally finite open covering of $S$ in $M$, each $s_\alpha$
is a smooth defining function for $S\cap U_\alpha$ in $U_\alpha$
with $ds_\alpha\neq 0$ in $U_\alpha$ and $\bar\partial_Ms_\alpha=0$
on $S\cap U_\alpha$, and moreover there exist
$\left\{g_{\alpha\beta}\right\}$, where $g_{\alpha\beta}$ is smooth
and nonzero on $U_\alpha\cap U_\beta$, such that
\form{s_\alpha\, = \, g_{\alpha\beta}s_\beta \qquad
\text{on}\quad U_\alpha\cap U_\beta\, .}
\smallskip
\noindent
{\sc Remark 1.}\quad Consistent systems of locally defining functions
always exist; in fact they arise naturally from a local
trivialization of the complex line bundle $N_SM$ and its
identification with a tubular neighborhood $V$, see (8.11).
\smallskip
\noindent
{\sc Remark 2.}\quad If our abstract $CR$ manifold $M$ is
$1$-pseudoconcave in a neighborhood of $S$, then any choice
of a system of local defining functions
$\left\{(U_\alpha,s_\alpha)\right\}$
with $s_\alpha$ satisfying (3.1) is automatically a consistent
system, according to Lemma 3.2.
\smallskip
Consider a form $\phi\in\Cal E^p(M\setminus S)$ and
$q=1,2,\hdots$.
\par
\noindent
{\sc Definition} (A) {\sl $\phi$ is said to have a pole of
(at most) order $q$ along $S$ iff: given any point $p\in S$,
there exists an open neighborhood $U$ of $p$ in $M$, and
a smooth local defining function $s$ of $S$ in $U$, with
$ds\neq 0$ in $U$ and $\bar\partial_Ms=0$ on $S\cap U$,
such that $s^q\phi\in\Cal E^p(U)$.}\par
(B) {\sl $\phi$ is said to be a {\bf semi-$CR$ meromorphic
form} with a pole of order (at most) $q$ along $S$
iff: it satisfies (A) with respect to some
consistent system of local defining functions.}
\smallskip
Note that (B) above implies that $\phi$ has
the local consistent representation in $U_\alpha$:
\form{\phi\, = \, \dfrac{\omega_\alpha}{s_\alpha^q}\, ,
\quad \omega_\alpha\in\Cal E^p(U_\alpha)\, .}
\par
We now come to the main point of this section: As we have seen,
given any $\phi\in\Cal E^p(M\setminus S,\Sigma)$ which is
closed, there exists another closed $\phi'\in
\Cal E^p(M\setminus S,\Sigma)$ cohomologous to $\phi$ which
has a simple pole along $S$. Then the residue class
$\Res[\phi]$ may be computed by
taking the class of $\res[\phi']$ in $\ho^{p-1}(S,\Sigma)$.
But the passage from $\phi$ to $\phi'$ by the route we have
(up to this point) developed is quite involved
and rather indirect. However when the closed $\phi$ is
semi-$CR$ meromorphic, and has a pole of finite order along $S$,
as in (B), we can do much better: It is possible to prescribe
an algorithm, which employs only elementary operations
on smooth differential forms, for the passage from $\phi$
to a cohomologous closed $\phi_1\in
\Cal E^p(M\setminus S,\Sigma)$ having a simple pole along $S$, and with
\form{\res[\phi_1] \in \Res[\phi]\, .}
\prop{Let $\phi\in\Cal E^p(M\setminus S,\Sigma)$ be a closed
semi-$CR$ meromorphic form having a pole of order $q\geq 2$
along $S$. Then we can construct semi-$CR$ meromorphic
forms $\hat\phi\in\Cal E^p(M\setminus S,\Sigma)$ and
$\rho\in\Cal E^{p-1}(M\setminus S,\Sigma)$, each having a
pole of order $(q-1)$ along $S$, such that
\form{\hat\phi\, = \, \phi\, - \, d\rho\qquad\text{and}\quad
\singsupp\left|_S\right.\hat\phi \subset
\singsupp\left|_S\right.\phi\, .}}
\dimo Differentiating (10.2) we obtain:
\form{0\, = \, \dfrac{d\omega_\alpha}{s_\alpha^q}
-q\dfrac{ds_\alpha\wedge\omega_\alpha}{s_\alpha^{q+1}}\quad
\text{in}\quad U_\alpha\setminus S\, ,}
from which we obtain
$ds_\alpha\wedge\omega_\alpha=0$ in $U_\alpha$ by continuity,
after wedging with $ds_\alpha$.
Hence by Cartan's lemma we can find $\theta_\alpha\in
\Cal E^p(U_\alpha,\Sigma)$ such that
$d\omega_\alpha=ds_\alpha\wedge\theta_{\alpha}$ in $U_\alpha$.
Using (10.5) we obtain
$ds_\alpha\wedge\left(q\omega_\alpha-s_\alpha\theta_\alpha\right)=0$
in $U_\alpha$; hence another application of Cartan's lemma
yields an $\eta_\alpha\in\Cal E^{p-1}(U_\alpha,\Sigma)$ such
that
$$q\omega_\alpha-s_\alpha\theta_\alpha=(q-1)ds_\alpha\wedge
\eta_\alpha\quad\text{in}\quad U_\alpha\, .$$
Therefore in $U_\alpha\setminus S$
the form $\phi$ can be represented
as
\form{\phi\, = \, \dfrac{\omega_\alpha}{s_\alpha^q}\, = \,
\frac{1}{q}\, \dfrac{\theta_\alpha}{s_\alpha^{q-1}}\, + \,
\frac{q-1}{q}\, \dfrac{ds_\alpha}{s_\alpha^q}\wedge
\eta_\alpha\, .}
 Let $\left\{\chi_\alpha\right\}$ be a smooth partition of
unity in a neighborhood of $S$, subordinated to
$\left\{U_\alpha\right\}$, and set $\chi=\sum_\alpha{\chi_\alpha}$.
Then we may write
$$\matrix\format\r\,&&\,\l\\
\phi&=&(1-\chi)\phi+\dsize\sum_\alpha{\chi_\alpha\phi}\\
&&\\
&=&(1-\chi)\phi+
\dsize\sum_{\alpha}{\chi_\alpha\dfrac{\omega_\alpha}{s_\alpha^q}}\\
&&\\
&=&(1-\chi)\phi+\frac{1}{q}\dsize\sum_\alpha{\left[
\dfrac{\chi_\alpha\theta_\alpha}{s_\alpha^{q-1}}
+(q-1)\dfrac{ds_\alpha}{s_\alpha^q}\wedge
\left(\chi_\alpha\eta_\alpha\right)\right]}\\
&&\\
&=&(1-\chi)\phi+\frac{1}{q}\dsize\sum_\alpha{\left[
\dfrac{\chi_\alpha\theta_\alpha}{s_\alpha^{q-1}}
-d\left(\dfrac{\chi_\alpha\eta_\alpha}{s_\alpha^{q-1}}\right)
+\dfrac{d\left(\chi_\alpha\eta_\alpha\right)}{s_\alpha^{q-1}}\right]}
\, .
\endmatrix$$
We get (10.4) by setting
$$\cases
\rho\, = \,
-\dsize\sum_{\alpha}{\dfrac{\chi_\alpha\eta_\alpha}{s_\alpha^{q-1}}}\, ,
\\
\hat\phi\, = \, (1-\chi)\phi+
\frac{1}{q}\dsize\sum_{\alpha}{\dfrac{\chi_\alpha\theta_\alpha+
d(\chi_\alpha\eta_\alpha)}{s_\alpha^{q-1}}}\, .
\endcases
$$
To see that $\rho$ ane $\hat\phi$ are semi-$CR$ meromorphic
with a pole of order $(q-1)$ along $S$, it suffices to observe that
$$
\cases
\rho\, = \, -\dfrac{\sum_{\beta}{g_{\alpha\beta}^{q-1}
\chi_\beta\eta_\beta}}{s_\alpha^{q-1}}\, ,\\
\\
\hat\phi\, = \, \dfrac{\sum_\beta{g_{\alpha\beta}^{q-1}
\left\{\chi_\beta\theta_\beta+d(\chi_\beta\eta_\beta)\right\}}}{q
s_\alpha^{q-1}}\, ,
\endcases
$$
in $U_\alpha\cap\{\chi=1\}$. The construction does not increase the
$\singsupp\left|_S\right.$; hence the proof is complete.
\thm{Let $\phi\in\Cal E^p(M\setminus S,\Sigma)$ be
a closed semi-$CR$ meromorphic form having a pole of order
$q\geq 2$ along $S$. Then we can construct a closed
semi-$CR$ meromorphic form $\phi^{(1)}\in\Cal E^p(M\setminus S,\Sigma)$
having a pole of order $q=1$ along $S$, such that
\form{\cases
\phi^{(1)}\, = \, \phi - d\rho\, ,\\
\res\left[\phi^{(1)}\right]\in\Res[\phi]\, ,\\
\supp \res[\phi^{(1)}]\subset
\singsupp\left|_S\right.\phi\, .
\endcases}}
\dimo Applying Proposition 10.1 $(q-1)$ times, we obtain
$$\matrix\format\c\,&&\,\l\\
\phi^{(q-1)}&=&\phi-d\rho^{(q-1)}\\
\phi^{(q-2)}&=&\phi^{(q-1)}-d\rho^{(q-2)}\\
\vdots&&\qquad\vdots\\
\phi^{(1)}&=&\phi^{(2)}-d\rho^{(1)}\, ,
\endmatrix
$$
where $\phi^{(j)}\in\Cal E^p(M\setminus S,\Sigma)$
and $\rho^{(j)}\in\Cal E^{p-1}(M\setminus S,\Sigma)$
are semi-$CR$ meromorphic and have poles of order $j$
along $S$. We obtain (10.7) with
$\rho=\rho^{(q-1)}+\rho^{(q-2)}+\cdots+\rho^{(1)}$,
proving the theorem.
\medskip
In $U_\alpha$ we have the consistent representations
\form{\phi^{(j)}=\dfrac{\omega_\alpha^{(j)}}{s_\alpha^j}\,,
\qquad \rho^{(j)}=\dfrac{\eta_\alpha^{(j)}}{s^j_\alpha}\, ,
\qquad 1\leq j\leq q-1\, ,}
with  $\omega^{(j)}_\alpha$, $\eta^{(j)}_\alpha$
smooth in $U_\alpha$. Thus in $U_\alpha$ we may write
\form{\dfrac{\omega_\alpha}{s_\alpha^q}\, = \,
\dfrac{\omega_\alpha^{(1)}}{s_\alpha}\, + \,
d\left\{\dfrac{\eta_\alpha^{(q-1)}}{s_\alpha^{q-1}}+
\dfrac{\eta_\alpha^{(q-2)}}{s_\alpha^{q-2}}+\cdots +
\dfrac{\eta_\alpha^{(1)}}{s_\alpha}\right\}\, ,}
and these local Laurent expansions are consistent
on $U_\alpha\cap U_\beta$.

\se{The Calculus of Residues}
In this section we assume that there is a neighborhood
$\Cal V$ of $S$
in $M$, in which $S$ has a global defining function $S$, with
$ds\neq 0$ in $\Cal V$ and $\bar\partial_Ms=0$ along $S$.
According to (9.4) we may, effectively, replace $M$ by $\Cal V$
when computing the residue class of a closed form $\phi$.
Hence there is no loss of generality if, in this section,
we take $M=\Cal V$. \par
We consider a semi-$CR$ meromorphic closed form $\phi\in
\Cal E^p(M\setminus S,\Sigma)$ with a pole of order $q$ along $S$,
having a global representation as
\form{\phi\, = \, \dfrac{\omega}{s^q}\, ,\qquad
\omega\in\Cal E^p(M,\Sigma)\, .}
In this situation we are able to convert the algorithm of the
previous section into a precise calculus of residues:
\par
consider the diffeomorphism:
\form{\CD
\nu: N_SM@>{\sim}>> V=M\, ,\endCD}
of the complex line bundle $N_SM$ with the tubular neighborhood.
It defines a foliation of $V$ with two-dimensional leaves
corresponding to the fibrs of $N_SM$. Let $L^*\subset T^*M$
be the set of covectors which annihilate
the vectors tangent to the leaves of the foliation. It is a rank
$(2n+k-2)$ real subbundle of $T^*M$. We denote by
$\Phi$ the $\Cal C^\infty$-subalgebra of the complexified
exterior algebra $\Cal E^*(M)$ generated by the smooth
sections of $L^*$ and by $d\bar s$. (Here $V=\Cal V=M$ is
chosen sufficiently small in order that $ds$ and $d\bar s$ are linearly
independent modulo $L^*$ at each point of $V$.) This gives us the
direct sum decomposition
\form{\Cal E^*(M)\, = \, \Phi\oplus ds\wedge\Phi\, .}
Using the decomposition (11.3) we define a $\Cal C^\infty(M)$-linear
operator
$\dfrac{d^0}{ds}$ on smooth forms:
\form{\dfrac{d^0}{ds}:\Cal E^p(M) \longrightarrow
\Phi\cap\Cal E^{p-1}(M)}
by saying that
\form{\dfrac{d^0\omega}{ds}\, = \, \beta}
iff $\omega\in\Cal E^p(M)$ decomposes as
$\omega=\alpha+ds\wedge\beta$ with
$\alpha\in\Phi\cap \Cal E^p(M)$ and
$\beta\in\Phi\cap \Cal E^{p-1}(M)$.
Note that
$\dfrac{d^0\omega}{ds}\left|_{\Sigma}\right.=0$
if $\omega\left|_{\Sigma}=0\right.$ because our tubular
neighborhood was adapted to $\Sigma$. We also define a $\Bbb C$-linear
operator
\form{\dfrac{d}{ds}:\Phi\cap\Cal E^p(M)
\longrightarrow\Phi\cap\Cal E^p(M)}
by saying that
$$\dfrac{d\lambda}{ds}=\beta$$
iff $\lambda\in\Phi\cap\Cal E^p(M)$ is such that
$d\lambda=\alpha+ds\wedge\beta$ with
$\alpha\in\Phi\cap\Cal E^{p+1}(M)$ and
$\beta\in\Phi\cap\Cal E^p(M)$. Again
$\dfrac{d\lambda}{ds}\left|_\Sigma\right.=0$
if $\lambda\left|_\Sigma\right.=0$.
\par
Note that if one choses local coordinates $(\Re s,\Im s,x)$
adapted to the foliation ($\{x=\roman{const}\}$ gives a leaf)
the operator (11.6) actually becomes the usual partial
differential operator
\form{\dfrac{\partial}{\partial s}\, = \,
\frac{1}{2}\left(\dfrac{\partial}{\partial\Re s}-\ssize{\sqrt{-1}}
\dfrac{\partial}{\partial\Im s}\right)}
acting on the coefficients of the form.
\par
We indicate the iterates of these operators by
\form{\cases
\dfrac{d^r}{ds^r}\, = \, \dfrac{d}{ds}
\hdots \dfrac{d}{ds}\, ,\quad\text{$r$ times}\, , \\
\dfrac{d^r}{ds^{r+1}}\, = \,
\dfrac{d^r}{ds^r}\dfrac{d^0}{ds}\, .
\endcases}
\prop{Let the closed form $\phi$ be as in (11.1) with $q\geq 1$.
Then
\form{\left.\frac{1}{(q-1)!}\dfrac{d^{q-1}\omega}{ds^q}\right| _S
\in\Res\left[\dfrac{\omega}{s^q}\right]\, .}}
\dimo Let
\form{\omega\, = \, \alpha+ds\wedge\beta}
where $\alpha\in\Phi\cap\Cal E^p(M,\Sigma)$ and
$\beta=\dfrac{d^0\omega}{ds}\in\Phi\cap\Cal E^{p-1}(M,\Sigma)$.
First we claim that if $q\geq 2$ then there are forms
$\tilde\alpha\in\Phi\cap\Cal E^{p}(M,\Sigma)$ and
$\eta\in\Cal E^{p-1}(M\setminus S,\Sigma)$ such that
\form{\phi+d\eta\, = \, \frac{1}{q-1}\dfrac{\tilde\alpha+
ds\wedge\frac{d\beta}{ds}}{s^{q-1}}\, .}
Indeed, since $d\phi=0$, we obtain
\form{q\,ds\wedge\alpha\, = \,
s\left(ds\wedge d\beta+d\alpha\right)\, .}
We decompose $d\alpha$ and $d\beta$ according to
(11.3), (11.6) as
\form{\matrix\format\r\,&&\,\l\\
d\alpha&=&\left.\left(d\alpha\right)\right| _\Phi+
ds\wedge\dfrac{d\alpha}{ds}\, ,\\
d\beta&=&\left.\left(d\beta\right)\right| _\Phi +
ds\wedge\dfrac{d\beta}{ds}\, ,
\endmatrix}
with $(d\alpha)\left| _\Phi\right.\in\Phi\cap\Cal E^{p+1}(M,\Sigma)$,
$\dfrac{d\beta}{ds}\in\Phi\cap\Cal E^{p-1}(M,\Sigma)$ and
$\dfrac{d\alpha}{ds},\, (d\beta)\left| _\Phi\right.
\in\Phi\cap\Cal E^p(M,\Sigma)$.
Then (11.12) gives
\form{\matrix\format\c\,&&\,\l\\
(d\alpha)\left| _\Phi\right.&=& 0\, ,\\
\alpha&=& \dfrac{s}{q}\left[
(d\beta)\left| _\Phi\right. + \dfrac{d\alpha}{ds}\right]\, ;
\endmatrix}
hence substituting in (11.10) we obtain
\form{\phi\, = \,
\dfrac{ds\wedge\beta}{s^q}+\frac{1}{q}
\dfrac{(d\beta)\left| _\Phi\right. +\frac{d\alpha}{ds}}{s^{q-1}}\, .}
From this we derive
\form{\phi+\frac{1}{q-1}d\left(\dfrac{\beta}{s^{q-1}}\right)\, = \,
\dfrac{\left(\frac{q}{q-1}(d\beta)\left| _\Phi\right.+
\frac{d\alpha}{ds}\right)+\frac{1}{q-1}ds\wedge\frac{d\beta}{ds}}{
s^{q-1}}\, .}
Thus we have proved our claim with $\tilde\alpha=
q(d\beta)\left| _\Phi\right. + (q-1)\dfrac{d\alpha}{ds}$.
\par
Next we claim that when $q=1$, then
\form{\res[\phi]\, = \, \left.\dfrac{d^0\omega}{ds}\right| _S\, .}
Indeed $d\phi=0$ gives (11.14) with $q=1$; hence (11.10),
(11.14) yields
\form{\phi\, = \, \dfrac{ds}{s}\wedge\beta+\left\{
\dfrac{d\alpha}{ds}+(d\beta)\left| _\Phi\right.\right\}\, ,}
which establishes our claim. We apply the
first claim $(q-1)$ times and the second claim once
to obtain the proposition.

\se{Iterated Residues}
We return to the general situation where $M$ is as in the beginning
of \S 3. But now we consider $m$ ($1\leq m\leq n$) polar
submanifolds $S_1$, $S_2$, $\hdots$, $S_m$ in $M$. We assume that
they are in {\it general $CR$ position} (i.e. have {\it normal crossing}).
This means that they are in general position and moreover that
their holomorphic tangent bundles $HS_1$, $HS_2$, $\hdots$, $HS_m$
are in general position in $HM$. Therefore the intersection of any
$r$ of them is a smooth closed $CR$ submanifold of type
$(n-r,k)$. In this section we set
\form{\cases
\check{S}&=S_1\cup S_2\cup\cdots\cup S_m\, ,\\
S&=S_1\cap S_2\cap\cdots\cap S_m\, ,\\
S_{(j)}&=\left(S_m\cap S_{m-1}\cap\cdots\cap S_{m-j+1}\right)
\setminus
\left(S_{m-j}\cup S_{m-j-1}\cup\cdots\cup S_1\right)\, ,
\endcases}
with $S=S_{(m)}$ and $M\setminus\check{S}=S_{(0)}$.
We consider also
$\Sigma=\Sigma_1\cup\Sigma_2\cup\cdots\cup\Sigma_\ell$, where each
$\Sigma_j$ is a smooth closed submanifold of $M$. We assume
that $S_1$, $S_2$, $\hdots$, $S_m$, $\Sigma_1$, $\Sigma_2$, $\hdots$,
$\Sigma_\ell$ are also in general position in the usual sense.
This implies, in particular, that the intersection of any subset
of the $\Sigma_1$, $\hdots$, $\Sigma_\ell$ is transversal to
the intersection of any subset of the $S_1$, $\hdots$, $S_m$.
\par
We define the $m$-th iterate
\form{\delta^m:\ho_{p-m}(S,\Sigma;\Bbb Z) \longrightarrow
\ho_p(M\setminus\check{S},\Sigma,\Bbb Z)}
of the coboundary homomorphism $\delta$ by
\form{\matrix\format\r&\r&\l&\l&\l&\l\\
\delta^m:\ho_{p-m}(S_{(m)},\Sigma;\Bbb Z) &@>\delta>>&
\ho_{m-p+1}(S_{(m-1)},\Sigma;\Bbb Z)&@>\delta>>&\cdots\\
&\cdots&@>>>\ho_{p-1}(S_{(1)},\Sigma;\Bbb Z)&@>\delta>>
&\ho_p(S_{(0)},\Sigma;\Bbb Z)\, .
\endmatrix}
Similarly we define the $m$-th iterate
\form{{\Res}^m:\ho^p(M\setminus\check{S},\Sigma)
\longrightarrow \ho^{p-m}(S,\Sigma)}
of the class residue homomorphism by
\form{\matrix\format\r&\r&\c&\l&\l\\
{\Res}^m:
\ho^p(S_{(0)},\Sigma)
&\overset{\Res}\to\longrightarrow&
\ho^{p-1}(S_{(1)},\Sigma)&
\overset{\Res}\to\longrightarrow&
\cdots\\
\hdots&
\longrightarrow
&\ho^{p-m+1}(S_{(m-1)},\Sigma)&
\overset{\Res}\to\longrightarrow&
\ho^{p-m}(S_{(m)},\Sigma)\, .
\endmatrix}
We may use Theorem 8.5 to obtain:
\thm{Let $h^*\in\ho^p(M\setminus\check{S},\Sigma)$ and
$\bar h\in\bar{\ho}_{p-m}(S,\Sigma;\Bbb Z)$. Then
\form{\dsize\int_{\overline{\delta^m h}}{h^*}
\, = \, \dsize\int_{\bar h}{\Res^m h^*}\, .}}
\medskip
\noindent
{\sc Remark 1.}\quad The homomorphisms $\delta^m$ ane
${\Res}^m$ defined here depend on the ordering
$S_1$, $S_2$,
$\hdots$, $S_m$. A permutation of the
$S_1$, $S_2$,
$\hdots$, $S_m$ affects $\delta^m$ and ${\Res}^m$ each by the
same $\pm$ sign, depending on the parity of the permutation.
\medskip
\noindent
{\sc Remark 2.}\quad Consider the situation where $m$ is maximal;
$m=n$. In this case $S$ is a {\it totally real} $k$-dimensional
$CR$ submanifold of type $(0,k)$ in $M$; moreover $S$ is
transversal to the
$CR$ structure of $M$, in the sense that
$T_xS\oplus H_xM=T_xM$ for every $x\in S$. Note that here
${\Res}^n$ could be regarded as the generalization to $CR$
manifolds of type $(n,k)$ of the Grothendieck residue (see [D]);
in which case we have $(k+1)$ kinds, as in (12.4) we may use
$\bar h\in \ho_j(S,\Sigma;\Bbb Z)$ for $j=0,1,\hdots,k$.

\se{The Calculus of Residues for Intersecting Poles}
We continue with the situation of \S 12 (polar submanifolds
$S_1$, $S_2$, $\hdots$, $S_m$ with normal crossings). But
in this section we assume that there is a neighborhood
$\Cal V$ of $S$ in $M$, in which each $S_j$ has a global
defining function $s_j$, with $\bar\partial_Ms_j=0$ along
$S_j\cap \Cal V$, for $j=1,\hdots ,m$. We may also assume that
the real and imaginary parts of
$ds_1$, $ds_2$,  $\hdots$,
$ds_m$ are linearly independent at each point of $\Cal V$.
As in \S 11 we can take $\Cal V=M$ and assume that $\Cal V=V$
is a tubular neighborhood of $S$, without any loss of generality.
\par
Having discussed iterated residues, we now turn our attention
to semi-$CR$ meromorphic forms having poles of finite order along
each $S_j$. Namely we consider the closed forms
$\phi\in\Cal E^p(M\setminus\check{S},\Sigma)$ having a global
representation as
\form{\phi\, = \,
\dfrac{\omega}{s_1^{q_1}s_2^{q_2}\cdots s_m^{q_m}}\, ,
\qquad\omega\in\Cal E^p(M,\Sigma)\, .}
Here we extend the calculus of residues of \S 11
to the {\sl multivariable} case. When $m>1$ we
write the linear operators in (11.5), (11.6), (11.8) as
$\dfrac{\partial^0}{\partial s_j}$,
$\dfrac{\partial}{\partial s_j}$,
$\dfrac{\partial^r}{\partial s_j^r}$,
$\dfrac{\partial^r}{\partial s_j^{r+1}}$
for each individual $S_j$. By composition we define
\form{\dfrac{\partial^{q_1+q_2+\cdots +q_m - m}\omega}{\partial
s_1^{q_1}\partial s_2^{q_2}\cdots\partial s_m^{q_m}}\, = \,
\dfrac{\partial^{q_1-1}}{\partial s_1^{q_1}}
\dfrac{\partial^{q_2-1}}{\partial s_2^{q_2}}
\hdots
\dfrac{\partial^{q_m-1}}{\partial s_m^{q_m}}
\omega\, , }
acting on $\omega\in\Cal E^p(M,\Sigma)$.
\prop{Let the closed form $\phi$ be as in (13.1)
with each $q_j\geq 1$. Then
\form{{\ssize{\frac{1}{(q_1-1)!(q_2-1)!\cdots (q_m-1)!}}}\left.
\dfrac{\partial^{q_1+q_2+\cdots +q_m-m}\omega}{\partial
s_1^{q_1}\partial s_2^{q_2}\cdots\partial s_m^{q_m}}\right|
_S \in {\Res}^m\left[
\dfrac{\omega}{s_1^{q_1}s_2^{q_2}\cdots s_m^{q_m}}\right]\, .}}
\dimo It suffices to observe that the operations described
in \S 11, with respect to a given $S_j$. commute with
the pullback to $S_i$, for $i\neq j$.\par
Note that the left-hand-side in (13.3)
can be computed in the usual sense of calculus.

\se{Abel's Global Residue Theorem}
In this section we present a generalization of the classical
theorem of Abel for a compact Riemann surface, along the lines
of Griffiths [G]. We return to the scenario at the beginning of
\S 12; so $S=S_1\cap S_2\cap\cdots\cap S_m$ and the
$S_1$, $S_2$, $\hdots$, $S_m$ are in general $CR$ position, etc.
\par
In particular $S$ is a $CR$ submanifold of $M$ of type $(n-m,k)$,
transversal to the $CR$ structure of $M$.
Its normal bundle $N_SM$ has the natural structure of a smooth
$\Bbb C$ vector bundle of rank $m$, via the identification
\form{N_SM\, =\, TM\left| _S\right.\left/ _{\dsize{TS}}\right.
\, \simeq\, HM\left| _S\right.\left/ _{\dsize{HS}}\right.\, .}
There are also natural maps $N_SM @>>> N_{S_j}M\left| _S\right.$
induced by projection onto the quotients, via the identifications
\form{\dfrac{N_SM}{N_{S}S_j}
\simeq
\dfrac{TM\left| _S\right.\left/ _{\dsize{TS}}\right.}{TS_j\left| _S
\right.\left/ _{\dsize{TS}}\right.}
\simeq
TM\left| _S\right.\left/ _{\dsize{TS_j}\left| _S\right.}\right.
\, = \,
N_{S_j}M\left| _S\right. \, .}
Therefore we get an isomorphism onto the Whitney sum
\form{N_SM \overset\sim\to\longrightarrow \bigoplus_{1\leq j\leq m}{
N_{S_j}M\left| _S\right.}\, .}
Analogous to \S 6 we construct a tubular neighborhood $W$ of $S$
in $M$ that is adapted to $\Sigma\cup\check{S}$:
\form{w:N_SM \overset\sim\to\longrightarrow W\, .}
Thus we have that
\form{w\left(\pi^{-1}\left([\Sigma\cup\check{S}]\cap S\right)\right)
\, = \, \left[\Sigma\cup\check{S}\right]\cap W\, ,}
where $\pi:W @>>> S$ denotes again the projection.
We fix a smooth Hermitian metric on the vector bundle
$N_SM$ such that the subbundles
$N_{S_1}M\left| _S\right.$,
$N_{S_2}M\left| _S\right.$, $\hdots$,
$N_{S_m}M\left| _S\right.$
are orthogonal, via the identification (14.3). A point $\zeta$
in the total space of $N_SM$ can be thought of
as $\zeta=(\zeta_1,\zeta_2\hdots ,\zeta_m)$, where
$\zeta_j$ is a point in $N_{S_j}M\left| _S\right.$ over
$\pi(\zeta)$. For the lenght of $\zeta$  at $x\in S$ we have
$|\zeta|_x^2=|\zeta_1|_x^2+|\zeta_2|^2_x+\cdots +
|\zeta_m|^2_x$.\par
Consider a relative $(p-m)$-cycle $\gamma$ in $(S,\Sigma)$.
As in \S 6 we identify $\gamma$ with a piecewise smooth map
\form{\hat\gamma:(P_\gamma,\partial P_\gamma) @>>> (S,\Sigma)\, ,}
where $P_\gamma$ is a finite polyhedron, of dimension
$(p-m)$, embedded in some Euclidean space $\Bbb R^N$.
Let $\gamma^*(N_SM)$ denote the pullback bundle over $P_\gamma$.
The map $\hat\gamma$ lifts to a smooth vector bundle
morphism
\form{\hat\gamma_{N_SM}:\hat\gamma^*(N_SM) @>>> N_SM\, ,}
giving by composition a map $f=w\circ\hat\gamma_{N_SM}$:
\form{f:\hat\gamma^*(N_SM) @>>> W\, .}
Next we consider the torus bundle $C^m_\gamma=C^m_\gamma(1)$
and the sphere bundle
$\Cal S^{2m-1}_\gamma=\Cal S^{2m-1}_\gamma(1)$ defined by
\form{C^m_\gamma\, = \,
\left\{ (y,\zeta)\in P_\gamma\times N_SM\,
\left| \, \pi(\zeta)=\gamma(y),\; |\zeta_j|_{\gamma(y)}=1,\;
j=1,\hdots,m\right.\right\}}
and
\form{\Cal S^{2m-1}_\gamma\, = \,
\left\{ (y,\zeta)\in P_\gamma\times N_SM\,
\left| \, \pi(\zeta)=\gamma(y),\; \sup_{1\leq j\leq m}{|\zeta_j|_{\gamma
(y)}}=1\right.\right\}\, .}
Note that
\form{\partial C^m_\gamma \, = \, C^m_\gamma \left| _{\partial P_\gamma}
\right.\, ,}
and
\form{\partial \Cal S^{2m-1}_\gamma\, = \,
\Cal S^{2m-1}_\gamma \left| _{\partial R_\gamma}\right. \, .}
Since our tubular neighborhood $W$ is adapted to $\Sigma$, we
may define
\form{\cases
\widehat{\delta^m\gamma}:(C^m_\gamma,\partial C^m_\gamma)
@>>> (M\setminus\check{S},\Sigma)\, ,\\
\widehat{\sigma^{2m-1}\gamma}:(\Cal S^{2m-1}_\gamma,
\partial\Cal S^{2m-1}_\gamma) @>>> (M\setminus S,\Sigma)\, ,
\endcases}
by $f$ restricted to $C^m_\gamma$ and $\Cal S^{2m-1}_\gamma$,
respectively.
Just as in \S 6 these define a smooth singular relative $p$-cycle
$\delta^m\gamma$ in $(M\setminus\check{S},\Sigma)$, and a smooth
singular relative $(p+m-1)$-cycle $\sigma^{2m-1}\gamma$ in
$(M\setminus S,\Sigma)$. Note that since $W$ was also adapted to
$\check{S}$, the $\delta^m\gamma$ just defined agrees with the
iterated $\delta^m\gamma$ in (12.2). \par
Set $U_j=M\setminus S_j$ for $j=1,2,\hdots ,m$; then $\Cal U=
\left\{U_j\right\}$ is an open covering of $M\setminus S$.
We shall associate to the covering $\Cal U$ two exact sequences,
one for {\v C}ech cochains of smooth differential forms, the other
for {\v C}ech chains of smooth singular chains.\par
First we denote by $\Cal C^{p,q}(\Cal U)$ the space of
alternating {\v C}ech $q$-cochains $g=\left(g_{i_0i_1...i_q}\right)$
with each $g_{i_0i_1...i_q}\in\Cal E^p(U_{i_0i_1...i_q},\Sigma)$,
where $U_{i_0i_1...i_q}=U_{i_0}\cap U_{i_1}\cap\cdots\cap U_{i_q}$.
The {\v C}ech coboundary operator $\check{\frak d}:\Cal C^{p,q}(\Cal U)
@>>>\Cal C^{p,q+1}(\Cal U)$ is given by
\form{
(\check{\frak d}g)_{i_0,i_1,...,i_{q+1}}=
\dsize\sum_{h=0}^{q+1}{
(-1)^h g_{i_0...\hat{i}_h...i_{q+1}}\left|
_{U_{i_0i_1...i_{q+1}}}\right.
}\, .}
Given $g\in\Cal E^p(M\setminus S,\Sigma)$ we let
$\epsilon^*g$ be the element of $\Cal C^{p,0}(\Cal U)$ defined by
\form{(\epsilon^*g)_i=g\left|_{U_i}\right.\, ,
\quad i=1,2,\hdots, m\, .}
Using a partition of unity subordinate to $\Cal U$, we obtain at once
the exactness of the sequence
\form{
0@>>>\Cal E^p(M\setminus S,\Sigma) @>{\epsilon^*}>>
\Cal C^{p,0}(\Cal U) @>{\check{\frak d}}>>
\Cal C^{p,1}(\Cal U) @>>>\cdots @>{\check{\frak d}}>>
\Cal C^{p,m-1}(\Cal U)@>>> 0\, .}
Second we denote by
$\Cal C_{p,q}(\Cal U)$ the $\Bbb Z$-module of alternating
{\v C}ech $q$-chains $\alpha=\left(\alpha_{i_0i_1...i_q}\right)$
with each $\alpha_{i_0i_1...i_q}$ being a smooth
singular $p$-chain in $(U_{i_0i_1...i_q},\Sigma)$.
The {\v C}ech boundary operator $\frak d:\Cal C_{p,q}(\Cal U)
@>>>\Cal C_{p,q-1}(\Cal U)$, $q\geq 1$, is defined by
\form{(\frak d\alpha)_{i_1i_2...i_q}=
\dsize\sum_{i_0=1}^m{\alpha_{i_0i_1i_2...i_q}}\, .}
Moreover we define $\epsilon_*:\Cal C_{p,0}(\Cal U) @>>>
\roman{Sing}_p(M\setminus S,\Sigma)$ by
\form{\epsilon_*\alpha=\dsize\sum_{i=1}^m{\alpha_i}\, .}
Here $\roman{Sing}_p(M\setminus S,\Sigma)$ denotes the space
of smooth singular relative $p$-chains in $(M\setminus S,\Sigma)$
with $\Bbb Z$-coefficients. By subdivision of the smooth
singular relative $p$-chains, we obtain the exactness of
the sequence
\form{0@>>>\Cal C_{p,m-1}(\Cal U)@>{\frak d}>>
\Cal C_{p,m-2}(\Cal U)@>>>\cdots
\cdots@>{\frak d}>>
\Cal C_{p,0}(\Cal U)@>{\epsilon_*}>>
\roman{Sing}_p(M\setminus S,\Sigma)@>>>0\, .}
\par
Next we introduce the duality pairing between $\Cal C^{p,q}(\Cal U)$
and $\Cal C_{p,q}(\Cal U)$:
\form{\dsize\int_{\alpha}{g}\, = \,
\dsize\sum_{1\leq i_0<i_1<\cdots<i_q\leq m}{
\dsize\int_{\alpha_{i_0i_1\hdots i_q}}{g_{i_0i_1\hdots i_q}}}\, ,}
for $g=\left(g_{i_0i_1\hdots i_q}\right)\in\Cal C^{p,q}(\Cal U)$
and
$\alpha=\left(\alpha_{i_0i_1\hdots i_q}\right)\in\Cal C_{p,q}(\Cal U)$.
The operators $\check{\frak d}$ and $\frak d$ are dual to
one another with respect to this pairing:
\form{\dsize\int_{\frak d\alpha}{g}\, = \,
\dsize\int_\alpha{\check{\frak d}g}\, ,\qquad
\text{for}\quad g\in\Cal C^{p,q}(\Cal U)\quad\text{and}\quad
\alpha\in\Cal C_{p,q+1}(\Cal U)\, .}
Likewise
\form{\dsize\int_{\epsilon_*\alpha}{g}\, = \,
\dsize\int_{\alpha}{\epsilon^*g}\, ,\qquad\text{for}
\quad g\in\Cal E^p(M\setminus S,\Sigma)\quad\text{and}\quad
\alpha\in\Cal C_{p,0}(\Cal U)\, .}
Now we use the exactness of (14.16) to show that there is
a homomorphism
\form{\sigma_*:\ho^p(M\setminus\check{S},\Sigma)
@>>> \ho^{p+m-1}(M\setminus S,\Sigma)\, .}
Indeed let $\phi$ be a closed form in $\Cal E^p(M\setminus\check{S},
\Sigma)$. We identify $\phi$ with an element $\phi^{(m-1)}
\in\Cal C^{p,m-1}(\Cal U)$ by
\form{\left(\phi^{(m-1)}\right)_{12\hdots m}\, = \, \phi\, ,}
as $M\setminus\check{S}=U_{12\hdots m}$. By the exactness of
(14.16) we can find $\phi^{(m-1)}$, $\phi^{(m-2)}$, $\hdots$,
$\phi^{(1)}$, $\phi^{(0)}$, with
$\phi^{(j)}\in \Cal C^{p,m-j-1}(\Cal U)$, and an
$\eta_\phi \in \Cal C^{p+m-1}(M\setminus S,\Sigma)$ with
$d\eta_\phi=0$ such that
\form{\cases
\phi^{(m-1)}\quad\text{is given by (14.23)}\\
\check{\frak d}\phi^{(m-1)}=\phi^{(m-1)}\\
\hdots\\
\check{\frak d}\phi^{(j)}=d\phi^{(j+1)}\quad\text{for}\quad
j=0,1,\hdots,m-3\\
\hdots\\
d\phi^{(0)}=\epsilon^*\eta_\phi\, .
\endcases}
Return to the consideration of our smooth singular relative
$(p-m)$-cycle $\gamma$ in $(S,\Sigma)$: for each
$0\leq j\leq m-1$ we define $\kappa_j\gamma\in\Cal C_{p+m-j-1,j}(\Cal U)$
by
\form{\left(\kappa_j\gamma\right)_{i_0i_1\hdots i_j}\, = \,
\left\{{\matrix\format\l\\
\text{the restriction of $f$ to}\quad \\
\quad\\
\left\{ (y,\zeta)\in\Cal S^{2m-1}_\gamma\,
\left| \, |\zeta_i|_{\pi(y)}=1\quad\text{if}\quad
i\in\{i_0,i_1,\hdots,i_j\}\right.\right\} \endmatrix}
\right.\,
,}
for $1\leq i_0<i_1<\cdots <i_j\leq m$. We note that
\form{\cases
\kappa_{m-1}\gamma=\widetilde{\delta^m\gamma}\leftrightarrow
\delta^m\gamma \\
\hdots	\\
\frak d\kappa_j\gamma=\partial\kappa_{j-1}\gamma\, ,
\quad\text{for}\quad 1\leq j\leq m-1\\
\hdots\\
\epsilon_*\kappa_0\gamma=\widetilde{\sigma^{2m-1}\gamma}
\leftrightarrow \sigma^{2m-1}\gamma\, .
\endcases}
Using Stokes' formula and (14.24), (14.26), (14.21), (14.22)
we obtain:
\form{\matrix
&\dsize\int_{\sigma^{2m-1}\gamma}{\eta_\phi}&=&
\dsize\int_{\epsilon_*\kappa_0\gamma}{\eta_\phi}&=&
\dsize\int_{\kappa_0\gamma}{\epsilon^*\eta_\phi}&=&
\dsize\int_{\kappa_0\gamma}{d\phi^{(0)}}\\
&&&\\
=&\dsize\int_{\partial\kappa_0\gamma}{\phi^{(0)}}&=&
\dsize\int_{\frak d\kappa_1\gamma}{\phi^{(0)}}&=&
\dsize\int_{\kappa_1\gamma}{\check{\frak d}\phi^{(0)}}&=&
\dsize\int_{\kappa_1\gamma}{d\phi^{(1)}} \\
&&\\
&\cdots&&\cdots&&\cdots&&\cdots\\
&&\\
=&\dsize\int_{\frak d\kappa_{m-1}\gamma}{\phi^{(m-2)}}&=&
\dsize\int_{\kappa_{m-1}\gamma}{\check{\frak d}\phi^{(m-2)}}&=&
\dsize\int_{\kappa_{m-1}\gamma}{\phi^{(m-1)}}&=&
\dsize\int_{\delta^m\gamma}{\phi}\, .
\endmatrix}
Thus we have established the following theorem,
which holds in the situation detailed in the
beginning of \S 12.
\thm{Let $\sigma_*:\ho^p(M\setminus\check{S},\Sigma)
@>>>\ho^{p+m-1}(M\setminus S,\Sigma)$ be the homomorphism
defined in (14.22). Then for any $h^*\in\ho^p(M\setminus\check{S},
\Sigma)$ and any class $[\gamma]\in\ho_{p-m}(S,\Sigma;\Bbb Z)$,
we have
\form{\left(2\pi\ssize{\sqrt{-1}}\right)^m
\dsize\int_{[\gamma]}{{\Res}^m h^*}\, = \,
\dsize\int_{[\delta^m\gamma]}{h^*}\, = \,
\dsize\int_{[\sigma^{2m-1}\gamma]}{\sigma_* h^*}\, .}}
\medskip
Now that we have Theorem 14.1, it is easy to derive a
generalization of Abel's theorem: We take $\Sigma=\emptyset$,
but otherwise $M$, $S$, $\check{S}$ are as in the beginning
of \S 12. When $M$ is compact, $S$ consists of a
finite number $S=Z_1+Z_2+\cdots +Z_N$ of connected
components $Z_i$, each being a smooth compact orientable
$CR$ submanifold of type $(n-m,k)$. A choice of orientation
of $M$ induces an orientation on each $Z_i$. This makes
each $Z_i$ into a smooth singular
$\left\{2(n-m)+k\right\}$-cycle in $M$.
\thm{Assume that $M$ is compact and $1\leq m\leq n$. If\par
\noindent
$\phi\in\Cal E^{2n+k-m}(M\setminus\check{S})$ is closed, then
\form{\dsize\sum_{i=1}^N{\dsize\int_{Z_i}{{\Res}^m[\phi]}}\, =
\, 0\, .}}
\dimo We consider $S$ with its natural orientation as a smooth
singular \par\noindent
$\left\{2(n-m)+k\right\}$-cycle. Then
$$\dsize\int_S{{\Res}^m[\phi]}\, = \,
\frac{1}{\left(2\pi\ssize{\sqrt{-1}}\right)^m}
\dsize\int_{\sigma^{2m-1}S}{\sigma_*[\phi]}\,=\, 0\, ,$$
because $-\sigma^{2m-1}S$ is the boundary of the complement
of a tubular neighborhood of $S$ in $M$, by Stokes' theorem.
\medskip
\noindent
{\sc Remark 1}. \quad For $k=0$ and $m=n$ we recover
a theorem of Griffiths [G]. The case $k=0$ and $m=n=1$
is the classical theorem of Abel for compact Riemann
surfaces.
\smallskip
\noindent
{\sc Remark 2.}\quad When the situation is as in \S 12 (each
$S_j$ has a global defining function in $V$), and $\phi$ is a
closed semi-$CR$ meromorphic form having a pole of finite
order along each $S_j$, as in (13.1), then (14.29) can be
written as
\form{\dsize\sum_{i=1}^N \left.
\dsize\int_{Z_i}{\dfrac{\partial^{q_1+q_2+\cdots +q_m-m}\omega}{\partial
s_1^{q_1}\partial s_2^{q_2}\cdots\partial s_m^{q_m}}}\right|
_{Z_i}\, = \, 0\, ,}
which is an expression involving "just calculus".
\se{Applications of the Abel's Theorem}
To illustrate the meaning of Theorem 14.2, we discuss in this
section only simple applications. We postpone any
investigation of complicated examples to a later date; it
would require indeed a lenghty discussion of the geometry
associated to different types of vector bundles over
$CR$ manifolds (see [HN3]). \par
Suppose $f_1,\, f_2,\, \hdots ,\, f_m$ are smooth
$CR$ functions on $M$, $1\leq m\leq n$, and consider
$S_j=\left\{x\in M\, \left| \, f_j(x)=0\right.\right\}$.
Set $S=S_1\cap S_2\cap \cdots\cap S_m$. We assume that there is
an open neighborhood $W$ of $S$ in $M$ in which each
$W\cap S_j$ is a smooth $CR$ submanifold of type $(n-1,k)$,
$df_j\neq 0$ in $W\cap S_j$. This means that each
$f_j$ has a simple zero along $W\cap S_j$, with
$df_j\wedge d\bar f_j\neq 0$ along $W\cap S_j$.
We assume also that the
$W\cap S_1$, $W\cap S_2$, $\hdots$, $W\cap S_m$ are in
general $CR$ position in $W$ (i.e. have normal crossings).
This implies that $S$ is a smooth $CR$ submanifold of type
$(n-m,k)$. Under the above assumptions we have:
\prop{Let $M$ be compact and oriented. Then we have the
period relations
\form{\dsize\sum_{i=1}^N\dsize\int_{Z_i}{\Theta}\, = 0\, ,}
for every smooth $\bar\partial_M$-closed form $\Theta$
of type $(n-m+k,n-m)$ on $M$.}
\dimo We consider the form
\form{g=\dfrac{df_1}{f_1}\wedge\dfrac{df_2}{f_2}\wedge
\dfrac{df_m}{f_m}\, ,}
defined on $M\setminus\check{S}$. Its ${\Res}^m$ is the
constant function $1$ on $S$. Our assumption on $\Theta$
implies that $\phi=g\wedge\Theta$ is closed.
${\Res}^m[\phi]=\left[\Theta\left| _S\right.\right]$.
Consider the class $[\eta_\phi]\in\ho^{2n+k-1}(M\setminus S)$
associated to $\phi$, where $\eta_\phi$ is as in (14.24).
Then we have by (14.27), (14.28) that
\form{\dsize\int_S{\Theta}\, = \,
\frac{1}{(2\pi\ssize{\sqrt{-1}})^m}\dsize\int_{\sigma^{2m-1}S}{
\eta_\phi}\, = \, 0\, ,}
because $\sigma^{2m-1}S$ is homologous to zero in
$M\setminus S$, and the proof is complete.
\medskip
\noindent
{\sc Remark.} Note that the compactness of $M$ does not prohibit
the existence of global $CR$ functions when $k>0$, as it would
when $k=0$. When $m=n$ we have that $S$ is a totally real $CR$
submanifold of type $(0,k)$, which is transversal to the $CR$
structure of $M$, and $\Theta$ is a $CR$ $k$-form on $M$.
\smallskip
Likewise under the above assumptions we have:
\prop{Let $M$ be compact and oriented. Then
\form{\dsize\sum_{i=1}^N\dsize\int_{Z_i}{
{\Res}^m\left[\dfrac{\omega}{f_1f_2\cdots f_m}\right]}\, = 0\, ,}
for every smooth closed form $\omega$ of type
$(n+k,n-m)$ on $M\setminus S$.}
\dimo
First observe that our hypothesis on $\omega$ implies that
$\phi=\omega\left/_{\dsize{f_1f_2\cdots f_m}}\right.$
is closed in $M\setminus\check{S}$. Consider the class
$[\eta_\phi]\in\ho^{2n+k-1}(M\setminus S)$ associated to
$\phi$ as in (14.24). Then we have as before that
\form{\dsize\int_S{{\Res}^m\left[\dfrac{\omega}{f_1f_2\cdots f_m}
\right]}\, = \, \dsize\int_{\sigma^{2m-1}S}{\eta_\phi}\, = 0\, ,}
since $\sigma^{2m-1}S\sim 0$ in $M\setminus S$.
\medskip
\noindent
{\sc Remark}. \quad When $m=n$ the $\omega$ above is a closed
section over $M\setminus S$ of the {\it canonical line bundle} of
$M$; i.e., a $CR$ $(n+k)$-form.

\enddocument